
\magnification=1200
\input amssym.def
\overfullrule0pt
\parindent=22pt
\def\EE{{\cal E}}
\def\D{{\Bbb D}}
\def\FF{{\cal F}}

\def\CC{{\cal C}}

\def\R{{\Bbb R}}
\def\N{{\Bbb N}}
\def\C{{\Bbb C}}

\def\11{{\bf 1\!\!}}
\def\ratop{\mathop{\hbox to .25in{\rightarrowfill}}\limits}

\def\nhang{\hangindent=4pc\hangafter=1}

\bigskip

\centerline{\bf Schwarz-type lemmas for solutions of 
$\bar\partial $-inequalities and } 

\centerline{\bf complete hyperbolicity of almost complex manifolds.}

\bigskip
\bigskip

\centerline{
by  Sergey Ivashkovich and Jean-Pierre Rosay%
\footnote*{Partly supported by NSF grant.}}

\footnote{} {AMS classification: 32Q60, 32Q65, 32Q45.}

\bigskip
\bigskip

\noindent {\bf 0. Introduction}

\medskip\noindent\sl 
0.1. Hyperbolic distance to a strictly pseudoconvex hypersurface.
\smallskip\rm 
In this paper an {\it almost complex manifold} $(X,J)$, means 
a smooth ($\CC^\infty$)
manifold $X$ with an almost complex structure $J$ of smoothness
at least $\CC^1$. When more smoothness is required 
by our proofs, it will be
specified in each statement. Through any point in any 
tangent direction,
there are local $J$-complex discs in $X$ (Theorem III in the pioneering 
paper [N-W], or see Appendix 1).
So, one can define the Kobayashi-Royden pseudo-norm of tangent vectors
and then the Kobayashi pseudo-distance on $(X,J)$.

The main question that will 
be studied in this paper is: for what domains $D\subset (X,J)$ is the almost
complex manifold $(D,J)$  complete hyperbolic (i.e. complete for the Kobayashi
distance in $D$)? 

Our first result is the following (see Proposition 2.1 for a local statement): 

\smallskip\noindent\bf
Theorem 1. {\it Let $D$ be a relatively compact strictly pseudoconvex domain of class
$\CC^2$ in an almost complex manifold $(X,J)$ and assume that $J$ is of class
$\CC^1$. Then either $(D,J)$ is complete hyperbolic or $D$ contains
a $J$-complex line.
}

\smallskip\rm Recall that a $J$-complex line in $D$ is an image of a 
non-constant $J$-holomorphic map $u:\C \to D$. 

There is one important case where the absence of $J$-complex lines in a
strictly pseudoconvex
$D$ is completely obvious. This is the case when $D$ possesses a global 
defining
strictly plurisubharmonic function i.e. there exists a ${\cal C}^2$-function 
$\rho $ in a neighborhood $V$ of $\bar D$ such that $D=\{ x\in V: \rho (x)<0\} $. 
Therefore we obtain

\smallskip\noindent\bf
Corollary 1. {\it Let $D$ and $(X,J)$ be as in Theorem 1 and suppose 
additionally that
$D$ possesses a global defining strictly plurisubharmonic function.
Then $(D,J)$ is complete
hyperbolic.
}

\smallskip\rm 

If $J$ is smooth enough and if the dimension of $X$ is 4 this result 
is due to Gaussier and Sukhov ([G-S]).
In higher dimensions Gaussier and Sukhov
prove the complete hyperbolicity of smoothly bounded strictly pseudoconvex
Stein domains with
almost complex structures sufficiently close to the initial Stein
complex structure. Their paper, which has been an inspiration, is based on rescaling 
and one then argues by contradiction. We directly deal with the needed estimates.

Corollary 1 implies the existence of bases of complete
hyperbolic neighborhoods of any point on $(X,J)$ provided $J\in \CC^1$.
Indeed, the  problem is local and therefore we can suppose that  
$X=\R^{2n}$ and $J(0)=J_{st}$. Then $\|\cdot \|^2$ is strongly 
$J$-plurisubharmonic near $0$. When the real dimension of $X$ is 4 the existence
of a basis of hyperbolically complete neighborhoods was proved in [D-I].

\smallskip It is worth to point out that strictly pseudoconvex domains in
almost complex manifolds are far different from 
those in complex ones, even if the almost complex structure is tamed
by some symplectic form. For example, in [McD] a symplectic 4-manifold 
$(X,\omega )$ together 
with a
relatively compact smoothly bounded domain $D\subset X$ is constructed such 
that
$\partial D$ is of contact type (and therefore is strictly pseudoconvex with respect 
to an appropriate
$\omega $-tamed $J$) but at the same time $\partial D$ is disconnected. In 
fact it has two
connected components. It is not clear however, whether this $D$ contains 
$J$-complex
lines or not.

\medskip\noindent\sl
0.2. Distance to a $J$-complex hypersurface. \rm

\smallskip A more precise analysis of McDuff's construction shows that this 
domain $D$ is defined as $\bar D=\{ x:t_1\le \rho (x)\le t_2\} $ where 
$\rho $ is some smooth function which is strictly pseudoconcave in $\bar D_1:=
\{ x:t_1\le \rho (x)\le t_0\} $ and is strictly pseudoconvex in 
$\bar D_2:=
\{ x:t_0\le \rho (x)\le t_2\} $ for some $t_0$. And $dd^c_j\rho =0$ on $\Gamma_
{t_0}$, where 
$\Gamma_{t_i}=\{ x:\rho (x)=t_i\}, i=1,2,3$. Therefore we obtain a domain, say 
$D_1$ with disconnected boundary $\Gamma_{t_1}\cup \Gamma_{t_0}$ such that 
one component is strictly pseudoconvex and another is Levi flat. Moreover, this 
$D_1$ doesn't contain $J$-complex lines.

In particular, for this reason we shall be interested in the sequel with hyperbolic distance
to Levi flat boundaries and to complex submanifolds. Therefore we now turn 
our attention to the hyperbolic distance to submanifolds in almost complex
manifolds. Let us state the problem more precisely.

Let $D$  be a
domain in an almost complex manifold $(X,J)$.
We do not assume that  $D$ is relatively compact nor that
$D$ is hyperbolic. 
A point $p\in \partial D$ is said to be at finite distance
from $q\in D$ if there is a sequence of points $q_j\in D$
converging to $p$ and whose Kobayashi distances to $q$ stay bounded.
Distances here are taken in Kobayashi pseudo-metric of $(D,J)$.
Otherwise we say that the distance is infinite.
Let $M$ be a closed submanifold of a domain $D$, of real codimension 1 or 2
(the case of higher codimension is trivial, see Remark at the end of \S5).  
For a point $p\in M$, we wish to investigate whether
there exist   points $q\in  D\setminus M$ at finite
Kobayashi distance from $p$ , in $D\setminus M$.
For codimension 2 our result is the following one:
\medskip

\noindent {\bf Theorem 2}. {\it 
Let $D$ be a be a hyperbolic domain in an almost complex manifold $(X,J)$, 
$J\in \CC^2$.
Let $M$ be a closed submanifold of $D$ of real codimension $2$ and of class
$\CC^3$. 

\itemitem{(2.A)} If $M$ is a $J$-complex hypersurface,then for every
$p\in M$ and
$q\in D\setminus M$ the Kobayashi distance from $q$ to $p$ is infinite.

 \itemitem{(2.B)} Conversely, if $p\in M$ and if the tangent space to
$M$ at $p$ is not $J$-complex, then for any neighborhood $D_1$ of
$p$ in $D$, there exists $p'\in D_1\cap M$ that is at finite distance from 
points in $D\setminus M$. 
\medskip
}
It is shown in section \S5 that at least if $J$ is of class $\CC^{2,\alpha}$,
instead of merely $\CC^2$, one can take $p'=p$, an immediate consequence of
Theorem 3.

Of course the existence of $J$-complex hypersurfaces is totally exceptional
unless the real dimension of $X$ is $4$.  However, Donaldson in [Do]
proved that every compact symplectic manifold admits symplectic
hypersurfaces in homology classes of sufficiently high degree, therefore giving
us almost complex structures with complex hypersurfaces on such manifolds.
Considering $J$-complex hypersurfaces comes naturally in the discussion of
real hypersurfaces.
Theorems 1 and 2 imply the following

\medskip\noindent {\bf Corollary 2}. {\it Let $J$ be a  $\CC^{2}$  
almost complex structure
in the neighborhood of the origin in $\R^{2n}$.
Suppose that there exists a $J$-complex $\CC^{3}$  hypersurface $M\ni 0$.
Then there exists a fundamental system of neighborhoods $\{ U_i\}$ of the
origin such that $(U_i,J)$ and $(U_i\setminus M,J)$ are complete hyperbolic in the sense
of Kobayashi.}
\medskip

In real dimension 4 this statement was proved in [D-I].

\medskip\noindent\sl
0.3. Distance to a Levi flat hypersurface.\rm
\smallskip
Now we turn to the case when $M$ is a real hypersurface. If $M$ is
Levi flat, i.e. foliated by $J$-complex hypersurfaces, then $D\setminus M$ is locally
hyperbolically complete as follows immediately from Theorem 2.
The inverse statement is not always true, in case of non integrable $J$.
In \S6 we construct the following:
\medskip
\noindent
{\bf Example.}  {\it There exists a real analytic almost complex 
structure $J$ on $\R^6$ such that $M=\R^5\times \{ 0\}$ is not Levi flat but the
Kobayashi pseudo-distance relative to
$(\R^6\setminus M,J)$ of any point in $\R^6\setminus M$ to $M$ is infinite.}

\medskip
Let a real hypersurface
$M\subset D$ be defined by $\rho =0$ (as usual $\nabla \rho \neq 0$ on $M$).  
A tangent vector $Y$ to $M$, at some point
$p$, is said to be {\it complex tangent\/} if $J(p)Y$ is also tangent to
$M$, at $p$.  Recall that $M$ is foliated by $J$-complex hypersurfaces
(i.e.\ $M$ has a Levi foliation) if and only if for every $p\in M$
and every $Y$ and $T$,
both complex tangent vectors to $M$ at $p$, $dd^c_J\rho (Y,T)=0$
(definitions will be recalled later).  

\smallskip\noindent\bf 
Theorem 3. {\it Let $D$ be a domain in an almost complex
manifold $(X,J)$. Assume that the closed real hypersurface $M\subset D$
is of class $\CC^2$, and J is of class $\CC^{3,\alpha}$
(for some $0<\alpha <1$).
If there exists a complex tangent vector $Y$ to $M$ at a
point $p$, such that $dd^c_J\rho (Y, JY)>0$, then that point $p$ is
at finite distance, in $D\setminus M$
from points in the region defined by $ \rho >0$.  
}

\medskip\rm
If $dd^c_J\rho (Y, JY)<0$, simply replace $\rho$ by $-\rho$. This theorem
implies that if a relatively compact smoothly bounded domain $D$
in almost complex manifold is complete hyperbolic then $D$ should be 
pseudoconvex (pseudoconvexity being defined by the Levi form).
\medskip

For non integrable almost complex structures the condition
$dd^c_J\rho (Y,JY)=0$ for all complex tangent vectors $Y$ does not
imply the (Frobenius) condition $dd^c_J \rho (Y,T)=0$ for
all complex tangent vectors $Y,T$. It is illustrated by the  above 
example. But there is a case when the above conditions are
equivalent. This is the case of dimension 4. We therefore have the following

\noindent {\bf Corollary 3.} {\it
Let $D$ be a
Kobayashi hyperbolic with respect to $J$ domain in $X$.
If $X$ has dimension 4, if the hypersurface $M\subset D$ is 
of class $\CC^3$ 
and $J$ is of class $\CC^{3,\alpha}$, then the following are equivalent:

(1) For every point $p \in M$ the Kobayashi distance, in $D\setminus M$,
from $p$ to any point in $D\setminus M$ is infinite.

(2) $M$ is Levi flat.}
\bigskip
\bigskip
The example in $\R^6$ mentioned above has the interesting feature that, although
$M=\R^5 \times \{ 0\}$ is not Levi flat, through any point in any
complex tangential direction, there is a $J$-complex curve lying entirely
in $M$. It would be interesting to know whether 
this geometric feature is actually equivalent to $dd^c_J\rho (Y,JY)=0$ for
all complex tangent vectors $Y$, and if both conditions imply complete hyperbolicity
of the complement.

\smallskip The structure of the paper is the following.

\smallskip\noindent
1. In \S1, we clarify some notations and recall basic facts from almost
complex geometry. It also contains some preliminaries, such as the existence
and basic properties of plurisubharmonic functions on almost complex manifolds
(including a quick proof
of well known basic results that also follow from standard elliptic theory).

\noindent
2. In \S2, we prove Theorem 1. The proof uses localization with the help of 
plurisubharmonic functions and then an appropriate Schwarz-type Lemma 2.3.

\smallskip\noindent
3. In \S3 we give an estimate of a Calderon-Zygmund integral and 
prove another Schwarz-type Lemma 3.2.

\smallskip\noindent
4. In \S 4 Lemma 3.2 is applied, and we prove Theorem 2.

\smallskip\noindent
5. Theorem 3 is proved in \S5.   
\smallskip\noindent
6. The example mentioned above is given in \S6. 
\smallskip\noindent
7. An appendix gathers some proofs and additional facts.

\bigskip

\noindent {\bf \S1. Some notations, definitions, and basic facts}.\rm

[McD-S] is a well known reference for almost complex manifolds.  An almost
complex manifold $(X,J)$ is a manifold $X$ of even real dimension $2n$,
with at each point $p$ an endomorphism $J=J(p)$ of the tangent space
satisfying $J^2= -\11$ . In this section, and in the next one we shall assume that
$J$ is of class $\CC^1$.

As usual $\CC^{k,\alpha}$ is used to denote spaces of maps whose
derivatives of order $\leq k$ are Holder continuous of order $\alpha$,
$k\in \N$, $0<\alpha<1$. $\CC^{k,\alpha}$ regularity of $J$
(i.e. of the map $p\mapsto J(p)$) is preserved by
$\CC^{k+1,\alpha}$ change of variables.

\bigskip \noindent
{\bf 1.a. $J$-holomorphic discs and the Kobayashi-Royden pseudo-norm.}
\bigskip
We shall study maps $u$ from an open set of $\C$ (always equipped
with the standard complex structure - so we will avoid the more
complete notation $\overline\partial_{J_{st},J}$) into $(X,J)$.  We set
$$
\overline\partial_Ju (Y)={1\over 2} [du (Y)+J(u)du (iY)],
$$
for any vector $Y$ tangent to $\C$ (at a point where $u$ is
defined).  The map $u\in \CC^0\cap L^{1,2}$ is $J$-holomorphic if
$\overline\partial_Ju =0$ a.e. If $J\in \CC^{k,\alpha }$ this implies that 
$u\in \CC^{k+1,\alpha }, k\ge 0$. In particular $\CC^1$-regularity
of $J$ implies $\CC^{1,\alpha }$-regularity of $u$, for any $\alpha\in (0,1)$. 
We shall see later that more is true: $u$ belongs then to some 
sub-Lipschitzian class
$\CC^{1,\phi}$ with $\phi (r)=r\ln{1\over r}$. It  will 
enable us to work under $\CC^1$-regularity of the structure.

$\D_R$ will denote the disc of radius $R$ in $\C$ centered at $0$, and
$\D =\D_1$ will be the open unit disc. Under our regularity assumption on 
$J$ through each point $p\in X$ in every direction $Y\in T_pX$ there
exists a $J$-complex curve, see \S 1.e. More precisely, there exists 
a $J$-holomorphic $u:\D_R \to X$ (for some $R>0$) such that $u(0)=p$ and 
$du(0)\left({\partial\over \partial x}\right)=Y$.
The Kobayashi-Royden pseudo-norm of the vector $Y\in T_pX$ on the
almost complex manifold $X$ is defined as $\| Y\|_K=\inf\Big\{{1\over
R}\mid \exists u:\D_R\to X$, $J$-holomorphic,
$u(0)=p,du(0)\left({\partial\over\partial x}\right) =Y\Big\}$. Another 
way to say the same is $||Y||_K=\inf\{ {1\over t}: \exists u:\D\to X, 
J-holomorphic , u(0)=p, du(0)\left({\partial\over \partial x}\right)=
tY\}$.

If $TX$ isequipped with some norm $\| \cdot \|$ then 
$\| Y\|_K=\inf\Big\{ {\|
Y\|\over \| du (0)\left({\partial\over\partial x}\right)\|}\mid
u:\D\to X,$ 

\noindent $~J-{\rm holomorphic},~u(0)=p~{\rm and}~du(0)({\partial\over 
\partial x})~ {\rm is ~parallel~ to }~ Y \Big\}$.

\smallskip The length $L$ of a path
$\gamma: [0~,~1]\to X$ is then defined by
$L=\int_0^1\| \dot\gamma (t)\|_k~dt$, where the integral is understood
as the upper integral, i.e. the
infimum of the integrals of the positive measurable majorants.

In case  the almost complex structure
$J$ is of class $\CC^1$ it is not clear that the function
$t\mapsto | \dot\gamma (t)\|_k$ is integrable or even measurable. 
But, at least if $J$ is of class $\CC^{1,\alpha}$, the Kobayashi-Royden
pseudo-norm is an upper semi continuous function on $TX$ (see Appendix 2), 
and the integral makes sense in the ordinary sense and is finite.

\bigskip
The pseudo-distance between two points is of course the 
infimum of the lengths of the paths joining these two points. If it is
a distance (separation property), the manifold is said to be
{\it hyperbolic}. Abusively we may say distance instead of
pseudo-distance.

\bigskip\noindent
{\bf 1.b. Complete hyperbolicity.}
We now state an elementary Lemma that will be used to
prove completeness for the Kobayashi metric.
It will be applied with $\delta (t)= Ct$ or
$Ct \log {1\over t}$ near 0 (C is a constant).
Think of $\chi$ in the Lemma as a complex coordinate function.
\medskip
\noindent {\bf Lemma 1.1.} {\it Let $D$ be a domain in an
almost complex manifold $(X,J)$. Let $p\in \partial D$.
Let  $\chi$ be either:
\smallskip
(a) a $\CC^1$ map from $\overline D$ into $\R^2$ with $\chi (p)=0$
and $\chi \neq 0$ on $D $, or
\smallskip
(b) a $\CC^1$ map for a neighborhood $U$ of $p$ into $\R^2$, such that
$\chi (p)=0$ and $\chi \neq 0$ on $U\cap \overline D \setminus \{ p\}$.
\smallskip
Let $\delta$ be a positive function defined on $(0~,+\infty)$
satisfying $\int_0^1 {dt\over \delta (t)} = +\infty$.
Assume that for every $J$-holomorphic map $u$ from $\D$ into
$D$, such that $u(0)$ is close to $p$
$$|\nabla (\chi \circ u)(0)|\leq \delta (|(\chi \circ u)(0)|).
\eqno (1.1)$$
Then, $p$ is at infinite Kobayashi distance from the points in $D$.}
\bigskip
\noindent
{\bf Proof.} We write first the proof for case (a). 
In the proof, the Kobayashi pseudo-norms are denoted by $|\cdot |_K$,
and $|\cdot |$ denotes Euclidean norm in $\R^2$.
\bigskip

Let $\gamma$ be a $\CC^1$ path in $D$ from a point 
$q_0$ (fixed) to a point $p_1$ (to thought of as close to $p$).
So $\gamma: [0~,~1]\to D$, $\gamma (0)=q_0$, $\gamma (1)=p_1$.
The Kobayashi length of $\gamma$ is 
$L=\int_0^1|\dot \gamma (t)|_K~dt$.
For $t\in [0~,~1]$, by definition of the Kobayashi metric, there exists 
a $J$-holomorphic
map $u_t:\D\to X$ with $u_t(0)=\gamma (t)$, and 

$$
{\partial u_t\over \partial x} (0)= 
{1 \over 2 |\dot \gamma (t)|_K}\dot \gamma (t)~.
$$
From this we get:

$$
{\partial (\chi\circ u_t)\over \partial x} (0)=d\chi_{\gamma (t)}({\partial u_t
\over \partial x}(0))={1\over 2|\dot\gamma(t)|_K}d\chi_{\gamma (t)}(\dot\gamma(t)
))=
{1\over  2 |\dot \gamma (t)|_K} {d\over dt}(\chi\circ \gamma) (t)~.
$$
Form (1.1) it follows that
$$|\dot \gamma (t)|_K \geq
{1\over 2} {|{d\over dt}(\chi\circ \gamma )(t)|
\over \delta(|(\chi\circ\gamma (t)|})
\geq {1\over 2} {-{d\over dt}(|(\chi\circ \gamma )(t)|)
\over \delta(|(\chi\circ\gamma (t)|)}~.$$
Finally, we get
$$L\geq {1\over 2} \int _{|\chi (p_1)|}^
                         {|\chi (q_0)|} {ds\over \delta (s)}~,$$
which is arbitrarily large if $p_1$ is close enough to $p$
(so $ |\chi (p_1)|\simeq 0$).
\bigskip
The proof of (b) follows from the above. Shrinking $U$ if needed, we can
assume that $\chi$ is defined on $\overline U\cap \overline D$ and
$\chi\neq 0$ on $\overline U\cap \overline D\setminus \{ p \}$.
We wish to estimate the length of a path from a point $q_0$ in $D$ to a point
$p_1$ as in the proof of (a). If this path is entirely in $U$, the proof of (a)
applies. Otherwise, let  $t_0$ the largest element in $[0,1]$
such that $\gamma (t_0)\notin U$ (assuming $U$ open). Then simply
apply the above estimates to the path from $\gamma (t_0)$ to
$p_1$ obtained by restricting $\gamma$ to $[t_0,1]$.
It is important to note that the non vanishing of $\chi$ on
$\overline U\cap \overline D\setminus \{ p \}$, and not only
on $\overline U\cap D$, insures that there exists $\epsilon >0$
such that $|\chi |\geq \epsilon$  on the boundary of $U$ (in $D$),
and therefore $|\chi (\gamma (t_0))|\geq \epsilon$ ($\epsilon$ not depending
on $\gamma$).

\smallskip\hfill{Q.E.D.}

\smallskip\noindent\bf
Remark 1. \rm  The Lemma covers the case of maps $\chi$ with values in $\R$,
such that $\chi |_{ D\cap U\setminus \{ p\}}<0$, and
$\chi (p)=0.$
\bigskip
\noindent {\bf Remark 2.}
In a hyperbolic manifold, complete hyperbolicity of a domain is
a purely local question. We have more:

An open subset $X_0$ in an almost complex manifold $(X,J)$ is called {\it
locally complete hyperbolic\/}  if for every $y\in \overline X_0$ there
exists a neighborhood $V_y\ni y$ such that $V_y\cap X_0$ is complete
hyperbolic.

Recall that an open subset $X_0$ of an almost complex manifold $X$ is
called {\it hyperbolically imbedded\/} into $X$ if for any two sequences
$\{ x_n\}$, $\{ y_n\}$ in $X_0$ converging to $x\in X$ and $y\in X$,
$x\neq y$,
respectively, one has that $\limsup_{n\to\infty} k_{X_0,J}(x_n, y_n)>0$.
Here $k_{X_0,J}$ denotes the Kobayashi pseudo-distance on the manifold
$(X_0,J)$.

It is worth observing that if $X_0$ is a relatively compact domain in $X$,
hyperbolically embedded into
$(X,J)$ and if $X_0$ is locally complete hyperbolic
then $(X_0,J)$ is complete hyperbolic, see [Ki].

The result is rather immediate (if $u$ maps the unit disc into $X_0$,
restrict to a
smaller disc and rescale), but it is very useful.  It reduces the problem
of completeness to a purely local problem (in hyperbolic manifolds).  We
will use it repeatedly without further explanations when switching to
local problems.

\noindent {\bf 1.c. Plurisubharmonic functions}.
\smallskip
If $\lambda$ is a function or vector valued map defined on $(X,J)$, 
$d^c_J\lambda$
is the $1$-form (vector valued) defined by
$$
d^c_J\lambda (Y)=-d\lambda (JY)
$$
for every tangent vector $Y$.  Notice that now the almost complex
structure is on the source space $X$.  Then $dd^c_J\lambda$ is defined by
usual differentiation.   
As usual, $\Delta$ will denote the Laplacian, 
($\Delta ={\partial^2\over\partial x ^2} + {\partial^2\over\partial y ^2}$).
The notation $d^c$ is relative to the standard complex structure on $\C$,
so for a function $h$ defined on an open set of $\C$: 
$d^ch
=-{\partial h\over\partial y} dx+ {\partial h\over\partial x} dy$.
 
In case of $\CC^2$-smoothness the following formula (1.2) follows immediately from
the chain rule and appeared in [De-2] and [Ha]. 

\medskip\noindent
{\bf Lemma 1.2.} {\it Let $J$ be a  $\CC^1$ almost complex structure
defined on an open set $\Omega \subset \R^{2n}$. Let $\lambda$ be a
$\CC^2$ function defined on $\Omega$. If $u:\D \to (\Omega ,J)$
is a $J$-holomorphic map, then:
$$\Delta (\lambda\circ u)~=~[dd^c_J\lambda ]_{u (z)}
\Big({\partial u\over \partial x}(z),[J(u(z)]({\partial u\over \partial x}(z)
\big)\Big). \eqno (1.2)$$
}
\bigskip
\noindent {\bf Comments}: Since $J$ is only of class $\CC^1$, $u$ need not
be $\CC^2$ (only $\CC^{1,\beta}$ for all $\beta<1$), so the left hand side
in the equation above has to be understood in the distributional sense.
But the right hand side is different. Since $dd^c_J$ involves only one derivative
of $J$, the right hand side makes sense pointwise (as suggested by
our insertion of $z$'s in the right hand side). The statement says that 
$\Delta (\lambda\circ u)$
is in fact (the distribution defined by) a continuous function,
although $\lambda\circ u$ need not be $\CC^2$.
\medskip
Although $u$ is not $\CC^2$, the Lemma should not be so
surprising. Consider the case of a genuine holomorphic change of variable
$\psi$, one has the formula
$\Delta (h\circ\psi)=\Delta h~|\nabla\psi|^2$, in which the second derivative
of $\psi$ plays no role.
\bigskip
\noindent {\bf Proof.} As already said $u$ is of class $\CC^1$,
so we can approximate $u$ by a sequence of smooth maps $v_j$
from $\D$ into $\Omega$, with $\CC^1$ convergence on compact sets.
\hfill\break
In the sense of currents on $\D$, $dd^c(\lambda\circ v_j)$ converges
to $dd^c(\lambda\circ u )= \Delta (\lambda\circ u)~dx\wedge dy$.
\hfill\break
We start by evaluating $d^c(\lambda\circ v_j )$. For any tangent vector
$T$ to $\D$, by definition of $d^c$:
$$d^c(\lambda\circ v_j ) (T) = - d(\lambda\circ v_j ) (iT)=
-d\lambda(d v_j(iT))$$
$$=~-d\lambda(Jdv_j(T))~+~\theta_j(T)~=~d^c_J\lambda(dv_j(T))~+~\theta_j(T),
\eqno (1.3)$$
with $\theta_j(T)=d\lambda\big( Jdv_j(T)-dv_j(iT)\big)$. 
Since $u$ is $J$-holomorphic
$Jdu(T)=du(iT)$, so $\theta_j$ tends uniformly to 0 on compact sets,
as $n\to\infty$. Consequently, $d\theta_j\to 0$ in the sense of currents.
Nothing more is to be used about $\theta_j$.
\medskip
Rewrite (1.3) $d^c(\lambda\circ v_j)=
v_j^\ast (d^c_J\lambda)+\theta_j$.
So
$$dd^c(\lambda\circ v_j)=v_j^\ast(dd^c_J\lambda)+d\theta_j.$$
Since $d\theta_j \to 0$, it follows that in the sense of currents on $\D$,
$$dd^c(\lambda\circ u )= {\rm lim}~dd^c(\lambda\circ v_j ) =
{\rm lim}~v_j^\ast(dd^c_J\lambda).$$
For the right hand side, the limit exists not only in the
sense of currents. Indeed
$v_j^\ast (dd^c_J\lambda)= h_j~dx\wedge dy$,
where 
$$h_j=v_j^\ast (dd^c_J\lambda) ({\partial\over \partial x},
{\partial\over \partial y})
= dd^c_J\lambda ({\partial v_j \over \partial x}, 
{\partial v_j \over \partial y}).$$
As $j\to \infty$, $h_j$ tends uniformly on compact sets to
$$dd^c_J\lambda ({\partial u \over \partial x}, 
{\partial u \over \partial y})~=~
dd^c_J\lambda ({\partial u \over \partial x}, 
J{\partial u \over \partial x})$$
(using again the $J$-holomorphicity of $u$).

The lemma is thus established.
\bigskip
\noindent {\bf Corollary 1.1.}  {\it Let $J$ be a $\CC^1$ almost complex structure
defined on an open set $\Omega \subset \C^n$. Let
$\lambda$ be a $\CC^2$ real valued function defined on $\Omega $. The following are
equivalent:

(1) For every tangent vector $Y$ to $\Omega$, $dd^c\lambda (Y,JY)\geq 0$.

(2) For every $J$-holomorphic map $u:\D \to \Omega$,
$\lambda\circ u$ is subharmonic.}
\bigskip
Let us repeat that since $dd^c_J$ involves only one derivative of $J$,
$dd^c_J\lambda$ makes clear sense. $(1)\Rightarrow (2)$ is an
immediate consequence of the Lemma. $(2)\Rightarrow (1)$ is also
an immediate consequence, taking into account that for every tangent vector
$Y$, at a point $q\in\Omega$, there exists a $J$-holomorphic map
$u: \D \to \Omega$ such that $u (0)=q$ and
${\partial u \over \partial x} (0)= t Y$, for some $t>0$.
\bigskip
If the equivalent conditions of the Corollary are satisfied by a
function $\lambda$, that function is said to be {\it $J$-plurisubharmonic}.
Of course (2) makes sense also for upper-semicontinuous $\lambda $ giving
us the general notion of a plurisubharmonic function. As usual, a function
is said to be strictly $J$-plurisubharmonic if locally any small
$\CC^2$ perturbation of that function is still $J$-plurisubharmonic.
\bigskip
There are two important examples of plurisubharmonic functions:

\medskip\noindent
{\bf Lemma 1.3.} {\it Let $(X,J)$ be an almost
complex manifold 
equipped with an arbitrary smooth Riemannian metric
with $J$ of class $\CC^1$.
Then for any $p\in X$, the function
$q\mapsto \big( {\rm dist}~ (q,p)\big)^2$ is strictly plurisubharmonic
near $p$.}
\medskip\noindent Proof. If $J$ is of class $\CC^{1,\alpha}$ the 
$J$-holomorphic discs are $\CC^{2,\alpha }$ and the result is trivial.
If $J$ is only of class $\CC^1$, it follows for (1.2) and from the
comments after Lemma 1.2.

The simple fact to be used is the following one: If $\mu$ is a continuous
function defined near 0 in $\R^2$ such that $\Delta \mu$  (in the sense of
distributions) is a continuous function, and such that for some $C>0$
$\mu(x,y)\geq \mu (0) + C(x^2+y^2)$ (near 0), then $\Delta \mu(0)>0$.
\bigskip

The next example less trivial is due to Chirka (not published).

\smallskip\noindent\bf
Lemma 1.4. {\it Let $J$ be an almost complex structure of class $\CC^1$ in the 
neighborhood of the origin in $\R^{2n}\simeq \C^n$. Suppose that $J(0)=J_{st}$. Then 
there exist a neighborhood $V\ni 0$ and a constant $A>0$ such that 
$\log |Z| + A|Z|$ is $J$-plurisubharmonic in $V$.
}\rm
\medskip
Consequently, for any point $p$ in a complex manifold $(X,J)$, with $J$ of class
$\CC^1$, there exists
a plurisubharmonic function $\lambda$ defined near $p$, continuous and finite
except at $p$ such that $\lambda (p)=-\infty$. Indeed we can take local coordinates
in which the almost complex structure coincides with the standard one at $p$.

\bigskip
\noindent{\bf Proof.} Using dilations the Lemma reduces to the following.
\smallskip

For $Z\in \R^{2n}~(=\C^n)$ set $u(Z)=|Z|+\log |Z|$.
We consider a continuously differentiable almost complex
structure $J$ defined on the unit ball $B$ in $\R^{2n}$.
We wish to prove that there exists $\epsilon >0$ such that if
$J(0)=J_{st}$ and $\|J-J_{st}\|_{\CC^1}\leq \epsilon$, then for
$Z\in B$, $Z\neq 0$,  and every tangent vector $Y$ at $Z$:
$$[dd^c_Ju(Z)](Y,JY)\geq 0~.$$
\bigskip
As previously $dd^c$ will be used for $dd^c_{J_{st}}$.
An elementary computation gives the following result:
$$dd^cu(Y,J_{st}Y)\geq {A\over |Z|} \|Y\|^2~,\eqno (1.4)$$
for some positive $A$ easy to determine. With complex notations:
it is equivalent to showing 
$\sum {\partial^2 u\over \partial z_j\partial\overline z_k} w_j\overline w_k
\geq {A\over 4 |Z|} |W|^2$.
\bigskip
Note that, using invariance under rotations, it is enough
to check the inequality at points $Z$ of the type
$Z=(z,0,\cdots ,0)$. the computations simplify and
one immediately gets a better non-isotropic estimate:
$$dd^cu(Y,J_{st}Y)\geq
 {A\over |z|} \|Y\|^2~+{A_1\over |z|^2}  \|Y'\|^2~,$$
with $Y=(Y_1,Y')\in \R^2\times \R^{2(n-1)}$.
\bigskip
Chirka's Lemma follows from a  simple perturbation argument in which we 
shall use that the first order (resp. second order) derivative of $u$
at $Z$ are of the order of magnitude of ${1\over |Z|}$ 
(resp. ${1\over |Z|^2}$).
We have:
$$dd^c_Ju~(Y,JY)~=~
 [d(d^c_J-d^c)u~(Y,JY)]~ 
+~[dd^cu~(Y,(J-J_{st})Y)]~+~ [dd^cu(Y,J_{st}Y)]~.$$
We now look at each of the terms on the right hand side.
For the last one $dd^cu(Y,J_{st}Y)$, we have the estimate $(1.4)$.
For the second one: 
$\|(J(Z)-J_{st})Y\|\leq \epsilon |Z|~\| Y\|$. So, using the estimates
on the derivatives of $u$, if $\epsilon$ is small enough:
$|dd^cu~(Y,(J-J_{st})Y)|\leq {A\over 4|Z|}\|Y\|^2$. 
\smallskip
We finally look at the first term.
Using $\R^{2n}$  coordinates $(t_1,\ldots ,t_{2n})$ for $Z\in B$, write
$$d(d^c_J-d^c)~=~\big( \sum \alpha_{j,k}^{l,m} 
{\partial^2\over \partial t_j \partial t_k}
+\sum \beta_j^{l,m}{\partial \over \partial t_j}\big)~ dx_l\wedge dx_m.$$
Since $J(0)-J_{st}=0$ and $\|J-J_{st}\|_{\CC^1}\leq \epsilon$, for some universal constant
$K$ we get $|\alpha_{j,k}^{l,m}|\leq K\epsilon |Z|$ 
and $| \beta_j^{l,m}|\leq K\epsilon$. 
Using the estimates on the derivatives of $u$, we again get that
if $\epsilon$ is small enough
$| d(d^c_J-d^c)u~(Y,JY)|\leq {A\over 4|Z|}\|Y\|^2$. Then one has:
$$dd^c_Ju~(Y,JY) \geq {A\over 2|Z|}\|Y\|^2~.$$
\hfill Q.E.D.

\bigskip

\bigskip\noindent{\bf
\S 1.d . Regularity of $J$-complex discs.}
\smallskip
\noindent  We want to make a few remarks about the regularity
of $J$-holomorphic discs since the proofs in the literature are not 
always pleasant. We shall 
always assume $u$ (as below) and $J$ to be of class $\CC^1$ at least, in which case
Formula (1.2) allows truly immediate proofs.

But it should be reminded that:
First of all if $J$ is merely continuous then $J$-holomorphic maps $u:\D\to X$ (a priori
they are $\CC^0\cap L^{1,2}$) are in $L^{1,p}$ for any $p$ and therefore in 
$\CC^{\alpha }$ for any $0<\alpha <1$ due to Sobolev embedding, see [IS-2]
Lemma 2.4.1. Then if $J$ is of class $\CC^\alpha$,  $u$ is of class
$\CC^{1,\alpha}$ (Theorem III in [N-W] and [Si]).
\smallskip\rm
\noindent {\bf  Remark 1.}
The proof of Lemma 1.2 used only the following:
$$
\left\{
\matrix{
\hbox{$\lambda$ is a $\CC^2$ function,}\hfill\cr
\hbox{ the almost complex structure $J$ is of class $\CC^1$,}
\hfill\cr
\hbox{ $u$ is a  $\CC^1$ $J$-holomorphic map.}\hfill\cr}
\right.\leqno(H)
$$
\medskip
 It may be
worth pointing out that
(using the simplest possible results on the regularity of the
standard Laplacian $\Delta$) 
the following is an immediate consequence of the validity of
formula (1.2) under the above hypotheses $(H)$ : \hfill\break
{\it If an almost complex structure $J$ on (some open set of) 
$\R^{2n}$ is of class $\CC^{k,\alpha}$, 
with $k\geq 1$ and $0<\alpha <1$, then any $\CC^1$ $J$-holomorphic 
map
$u : \D \to  (\R^{2n},J)$ is $\CC^{k+1,\alpha}$
regular. If $J$ is only $\CC^1$, $u$ is of class
$\CC^{1+\beta}$ for all $\beta<1$}. More is true (see Remark 3)
and will be needed later.
\hfill\break
Here, we sketch the argument:   
\medskip
1) Assume that $J$ is of class $\CC^1$. Let $\lambda$ be any
smooth function defined on $\R^{2n}$.

By $(1.2)$, $\Delta (\lambda \circ u)$ is a locally bounded function.
It immediately follows that $\lambda\circ u$ is of class $\CC^{1,\beta}$
for all $0<\beta <1$. Taking $\lambda$ to be the coordinate functions,
one sees that $\lambda$ itself is of class $\CC^{1,\beta}$.

2) Assume now that $J$ is of class  $\CC^{k,\alpha}$ ($k\geq 1~,~
0<\alpha <1$). By 1) we already know that $u$ is
of class $\CC^{1,\alpha}$.
Then, by repeating the argument, one sees that if $u$ 
is of class $\CC^{r,\alpha}$ with $1\leq r \leq k$, then 
$u$ is of class $\CC^{r+1,\alpha}$.
\medskip\noindent
Indeed, if $\lambda$ is any smooth function on $\R^{2n}$,
(1.2) shows that $\Delta (\lambda\circ u)$ is of class
$\CC^{r-1,\alpha}$. So by regularity of the Laplacian
$\lambda\circ u$ is $\CC^{r+1,\alpha}$.

\bigskip\noindent
{\bf Remark 2.}
The other result of basic elliptic theory that
is used several times is that if on a bounded open set $\Omega$
in $\R^{2n}$, an almost complex structure $J$ is close enough 
(depending on $\beta$ below, $0<\beta<1$) to
$J_{st}$ in $\CC^1$ norm (resp. close enough in $\CC^{k,\alpha}$ norm, with
$k\geq 1$), then the $J$-holomorphic maps $u$ from
$\D$ into $\Omega$ have uniform
$\CC^{1,\beta}$ (resp. $\CC^{k+1,\alpha}$) bounds on smaller discs.
It is at the root of local hyperbolicity, see Lemma 3.3.
We shall restrict our discussion to almost complex structures which are at
least of class $\CC^1$ (see [Si], for lower regularity).

\bigskip
We wish to mention that formula (1.2)
can also be used for proving the above result.
One first needs an initial regularity  result giving $L^p$ bounds
for $\nabla u$, to be used with $p>4$.  Then, proceeding as above,
a first application of (2.2)
gives $L^r$ bounds for $\Delta u$ with $r={p\over 2}>2$, therefore
H\" older
$(1-{2\over r})$ estimates for $\nabla u$, i.e.  
$\CC^{1,1-{2\over r}}$ bounds for $u$. After that, (1.2)
gives $\CC^{k+1,\alpha}$ bounds if $J$ is of class $\CC^{k,\alpha}$
and if $\CC^{k,\alpha}$ bounds are already known for $u$. 
\bigskip
It happens that $L^p$ bounds for $\nabla u$, for arbitrary
$1< p < \infty$, are extremely easy to get.  Provided that $\epsilon$
is small enough depending on $p$ (for $p=2$, $\epsilon <1$), 
they follow from the simple
differential inequality:
$$\left| {\partial u \over \partial \overline z}\right|
\leq
\epsilon \left| {\partial u \over \partial  z}\right|.$$
The argument is well known. Let 
$\chi \in \CC^\infty_0(\D)$, satisfying $0\le \chi \le 1$. We have:
$$\left| {\partial \chi u \over \partial \overline z} \right|
\leq 
\left| \chi {\partial u \over \partial \overline z} \right| +
\left| {\partial \chi  \over \partial \overline z} u \right|
\leq 
\epsilon \left| {\partial \chi u \over \partial z} -
{\partial \chi \over \partial z} u \right|
+ \left| {\partial \chi  \over \partial \overline z} u \right|
\leq
\epsilon \left| {\partial \chi u \over \partial z} \right|
+(1+\epsilon) K\vert u\vert ,$$
with $K= {\rm Sup}~|\nabla \chi| .$
Since $ {\partial \chi u \over \partial z} 
={-1\over \pi z^2}\ast {\partial \chi u \over \partial \overline z}$,
the theory of singular integral gives
$\| {\partial \chi u \over \partial z}\|_{L^p}
\leq C_p \| {\partial \chi u \over \partial \overline z}\|_{L^p}
$. So,  one has
$$\| {\partial \chi u \over \partial \overline z}\|_{L^p}
\leq
\epsilon C_p
\| {\partial \chi u \over \partial \overline z}\|_{L^p}
+(1+\epsilon ) K\| u\|_{L^p}.$$
If $\epsilon$ is chosen small enough so that $\epsilon C_p<1$,
we get:
$$\| {\partial \chi u \over \partial \overline z}\|_{L^p} \leq
{(1+\epsilon ) K\over 1-\epsilon C_p} \| u\|_{L^p}.$$
Then, $$\| {\partial \chi u \over \partial z} \|_{L^p} \leq C_p
{(1+\epsilon ) K\over 1-\epsilon C_p} \| u\|_{L^p}.$$
Finally,
$$\| \nabla \chi u\|_{L^p} \leq (1+C_p) {(1+\epsilon ) K\over 1-\epsilon C_p} \| u\|_{L^p}.
$$
\bigskip
We will need the following that is easily obtained
from Remark 2 above.
\bigskip
\noindent {\bf Lemma 1.5 .} {\it Let $\Omega$ be an open subset 
of $(\R^{2n},J)$, $J$ of class $\CC^1$. 
Let $K$ be a compact subset of $\Omega$. There exists $\delta>0$,
such that: for every $r\in [0,1)$
there exists $C>0$
such that if $u:\D \to \Omega$ is a $J$-holomorphic disc
with $u(\D)\subset K$, then
$$
|\nabla u(z)| \leq C {\rm sup}_{|t|<1}|u(t)-u(0)|,\eqno(1.5)
$$
if ${\rm sup}_{|t|<1}|u(t)-u(0)|\leq \delta$ and
$|z|\leq r$.}
\bigskip
\noindent Proof. Depending on $u(0)$, one can make a linear change of variables
such that in the new coordinates $J(u(0))=J_{st}$.
One can choose the linear maps for changing variables so that
their norms and the norm of their inverses are uniformly bounded
for $u(0)\in K$.
\smallskip
After such a change of variables, set $M={\rm Sup}~|u(z)-u(0)|$
($z\in \C$, $|z|<1$). Set $u_M(z)={u(z)-u(0)\over M}$.
Then $u_M$ is a $J_M$-holomorphic map from $\D$ into the unit ball
in $\R^{2n}$, for the almost complex structure $J_M$ that is the
push-forward of $J$ under the map $Z\mapsto {1\over M}(Z-u(0))$.
If $M$ is small enough, $J_M$ is close to $J$ in the $\CC^1$ sense
and Remark 2 applies. So, for $|z|\leq r$ one gets $|\nabla u_M(z)|\leq C$
for some absolute constant $C$, hence $|\nabla u(z)|\leq CM$, as desired.
\bigskip
\noindent {\bf Remark 3.}
Let $\CC^{1,\phi}$ (see \S3 for an explication of the notation)
be the  class of 
continuously differentiable functions or maps that locally satisfy the
following sub-Lipschitzian condition
$$|\nabla f(z')-\nabla f(z)|\leq C |z'-z|~\log{1\over |z'-z|}~,$$
for some constant $C$. 
So, the gradients of functions in  $\CC^{1,\phi}$ are 
better than H\" older, but not quite Lipschitzian.
\bigskip\noindent
{\bf Lemma 1.6.}  {\it Any $J$-holomorphic disc (that in our proof we assume to be
$\CC^1$) is in the class  $\CC^{1,\phi}$ if $J$ is $\CC^1$.}
\medskip Lemma 1.6 follows immediately from
the observations in 1) in  Remark 1 and from the elementary Lemma:
\bigskip
\noindent {\bf Lemma 1.7.} {\it Let $g$ be a function defined on some
open set in $\R^2\simeq \C$. If $\Delta g$ (in the sense of distributions)
is a bounded function, then $g\in \CC^{1,\phi}$.}
\bigskip
\noindent {\bf Proof.} Assume that $g$ is defined near 0. Let $\chi$
be a smooth cut off function such that $\chi \equiv 1$ near 0.
We have $g\chi=\Delta (g\chi )\ast{1\over 2 \pi} \log|z|$.
So, near 0:
$${\partial g\over \partial z}=\Delta (g\chi )\ast {1\over \pi z}~{\rm and}~
{\partial g\over \partial \overline z}=\Delta (g\chi )\ast {1\over \pi \overline z}.$$
We concentrate on ${\partial g\over \partial z}$, the case of
${\partial g\over \partial \overline z}$ being similar.
\medskip
Write $\Delta (g\chi )$ as the sum of a bounded function $v$ with compact support
and of a distribution $T$ with support away from 0. Then, $T\ast {1\over \pi z}$
is $\CC^\infty$ near 0. So we only need to show an estimate of the type:
$$|(v\ast {1\over \pi z} )(z_0 +t)-(v\ast {1\over \pi z} )(z_0)|
\leq C |t|\log{1\over |t|}.$$
Checking that estimate is immediate. Write
$$(v\ast {1\over \pi z} )(z_0 +t)-(v\ast {1\over \pi z} )(z_0)~=$$
$$\int {v(z_0-\zeta)\over \zeta + t}~dxdy(\zeta ) -
\int {v(z_0-\zeta)\over \zeta}~dxdy(\zeta ).$$ 
Then, simply use the estimates:
$$\int_{|\zeta |<2 |t|} {1\over |\zeta |}~{\rm and}~
\int_ {|\zeta |<2 |t|} {1\over |\zeta +t|}~=~O(|t|),$$
and
$$\int_{|\zeta |>2 |t|} |{1\over \zeta +t}-{1\over \zeta}|~\leq~
2 |t| \int_{|\zeta |>2 |t|} {1\over |\zeta|^2}~=~O(|t|\log |t|)~.$$
\hfill{Q.E.D.}
\bigskip
\noindent{\bf 1.e. Jets of $J$-holomorphic discs.}
\bigskip
The case of 1-jets ($k=1$)is simply the case of discs with
prescribed tangent.
In the proof of Theorem 3, we shall also need the case of 2-jets,
with dependence on parameters. But we 
state the general case of jets of arbitrary order.
\bigskip
\noindent {\bf Proposition 1.1.} {\it Let $k\in \N$, $k\geq 1$, and
$0<\alpha<1$. Let $J$ be a $\CC^{k-1,\alpha}$ almost complex structure 
defined near 0 in $R^{2n}$. For any $p\in \R^{2n}$ close enough to 0,
and every $V=(v_1,\ldots ,v_k)\in (\R^{2n} )^k$ small enough,
there exists a $\CC^{k,\alpha}$ $J$-holomorphic map $u_{p,V}$ from $\D$ 
into $R^{2n}$ such that $u_{p,V}(0)=p$, and ${\partial^l u_{p,V} \over 
\partial x^l}(0)=v_l$, for any $1\leq l\leq k$.
If the structure $J$ is of class $\CC^{k,\alpha}$, then
$u_{p,V}$ can be chosen with $\CC^1$ dependence (in $\CC^{k,\alpha}$)
on the parameters $(p,V)$ in $\R^{2n}\times (\R^{2n} )^k.$
}\bigskip
Here we will give the proof assuming $\CC^{k,\alpha}$ regularity of
$J$. It is a rather simple consequence of the implicit function theorem. 
The proof of the existence of discs under $\CC^{k-1,\alpha}$ regularity requires 
a trick found by [N-W] for one-jets, and the Schauder 
fixed point Theorem. It is given in the Appendix. We do not know whether continuous 
dependence on parameters can then be achieved (there are possibly related
examples indicating that it is conceivable that it fails).
\bigskip
Before the proof of the Proposition, we start with some preliminaries.
\bigskip
\noindent {\bf 1.e.1.} Re-writing of $J$-holomorphicity. 
\bigskip
On $\R^{2n}$, we consider an almost complex
structure $J$ and the standard almost complex structure $J_{st}$
(corresponding to multiplication by $i$ in the identification
of $\R^{2n}$ with $\C^n$). By definition
$$\overline\partial_J u ({\partial\over \partial x})~=~
{\partial u\over \partial x}
+J(u){\partial u\over \partial y}.\eqno (1.6)$$
With some abuse of notations ($\R^{2n}$ notations on the left
hand side and $\C^n$ notations on the right hand side)
$${\partial u \over \partial x}
= {\partial u \over \partial z}+{\partial u \over \partial \overline z}
~;~
{\partial u \over \partial y}= {1\over i} (-{\partial u \over \partial z}
+ {\partial u \over \partial \overline z}) =
-J_{st}  (-{\partial u \over \partial z}
+ {\partial u \over \partial \overline z}).$$
By multiplication on the left by $J(u)$, $(1.6)$ gives:
$$J(u) \overline\partial_J u ({\partial\over \partial x})=
[J(u)+J_{st}]{\partial u\over \partial \overline z}
+[J(u)-J_{st}]{\partial u\over \partial z}~.$$
We shall restrict our attention to almost complex structures
$J$ such that at each $Z\in X$, $J(Z)+J_{st}$ is invertible,
which happens in particular if $J(Z)-J_{st}$ has operator norm
$<1$.\hfill\break
Then, set:
$$Q_J(u)=[J(u)+J_{st}]^{-1}[J(u)-J_{st}]~.$$
The equation for $J$ holomorphicity of $u$ becomes
$${\partial u\over \partial \overline z} 
+Q_J(u){\partial u \over \partial z}~=~0~.\eqno (1.7)$$
Indeed since $\overline \partial_Ju$ is simply the $\C-J$ anti-linear
part of $du$, $\overline \partial_Ju=0$ if and only if
$\overline \partial_Ju ({\partial \over \partial x})=0.$
\bigskip
\noindent {\bf 1.e.2.} The Cauchy-Green operator $T_{CG}$ and the operator
$\Phi_J$.\bigskip
For a complex valued function $g$ or a map  $g$ with values in a complex
vector space $g$ continuous on $\overline \D$, and $z\in \C$ with
$|z|\leq 1$, we set:
$$T_{CG}(g)(z)~=~(g\ast {1\over \pi \zeta})(z)~=~
{1\over \pi} \int_D{g(\zeta)\over z-\zeta}~dxdy(\zeta )~.$$
We shall need the classical properties of $T_{CG}$ (see Appendix 4):
\smallskip
(a) If $g\in \CC^{k,\alpha}(\overline \D ),~k\in \N~,0<\alpha<1~,$
then $T_{CG}g\in \CC^{k+1,\alpha}(\overline \D )$.
\smallskip
(b) ${\partial \over \partial\overline z}[T_{CG}(g)]=g$. (on $\overline \D$.)
\bigskip
Let $k$ be an integer $\geq 1$. 
Assume that $J$ is a $\CC^{k,\alpha}$ almost complex structure 
defined on $\R^{2n}$, such that $J(z)+J_{st}$ is invertible for all
$z\in \R^{2n}$.
We define the operator $\Phi_J$ from
$\CC^{k,\alpha}(\overline \D,\R^{2n})$ into itself by:

$$
\Phi_J(u)~=~({\bf 1}-T_{CG}Q_J(u){\partial \over \partial z})u~.
\eqno(1.8)
$$
$\Phi_J$ is a continuously differentiable map
from $\CC^{k,\alpha}(\overline \D,\R^{2n})$ into itself, whose 
derivative at the point $u\in \CC^{k,\alpha}(\overline \D,\R^{2n})$
is the map:

$$\delta u \mapsto
\Phi_J(\delta u)-T_{CG}S(\delta u) {\partial u\over \partial z},$$
where $S$ is the operator defined by differentiation of
$u\mapsto Q_J(u)$ (as a map from $\CC^{k,\alpha}$ into
$\CC^{k-1,\alpha}$ - $T_{CG}$ re-gaining one derivative), 
i.e. with obvious notations:
$$S(\delta u)
= [J(u)+J_{st}]^{-1} DJ(\delta u)~-~
[J(u)+J_{st}]^{-1} DJ(\delta u) [J(u)+J_{st}]^{-1}[J(u)-J_{st}].$$
It is for this differentiation that it is not enough
that $J$ be of class $\CC^{k-1,\alpha}$.
If $J=J_{st}$, $\Phi_J$ is the identity mapping.
On any fixed ball in $\CC^{k,\alpha}$,
$\Phi_J$ is a small (non-linear) perturbation of the identity if
$J$ is close to $J_{st}$ in $\CC^{k,\alpha}$ topology. Finally note that
equation $(1.7)$ and (b) show that $u$ is $J$-holomorphic if and only if
$\Phi_J(u)$ is holomorphic in the ordinary sense.
\bigskip
\noindent {\bf 1.e.3} Proof of Proposition 1.1, assuming $\CC^{k,\alpha}$
smoothness of $J$. 
\bigskip
Of course, we can assume that $J(0)=J_{st}$. After linear change of
variables (dilations), and cut off of $J-J_{st}$, we can assume that 
$J$ is defined on $\R^{2n}$, and as close as we wish to $J_{st}$
in $\CC^{k,\alpha}$ topology. 
For $(q,W)\in \R^{2n}\times (\R^{2n} )^k$
($W=(w_1,\ldots ,w_k)$), define the map $h_{q,W}$
on $\overline D$ by
$$ h_{q,W}(z)=q+\sum_{l=1}^k{1\over l!}z^lw_l~.$$
Fix $R>0$ so that for any $(q,W)$ in the ball ${\bf B}_R$ of radius $R$
in $\R^{2n}\times (\R^{2n} )^k$, $h_{q,W}$ is in the unit ball of
$\CC^{k,\alpha}(\overline \D,\R^{2n})$.
If $J$ is close enough to $J_{st}$, $\Phi_J^{-1}$ is defined
as a map from the unit ball of $\CC^{k,\alpha}(\overline \D,\R^{2n})$ 
into $\CC^{k,\alpha}(\overline \D,\R^{2n})$. For $(q,W)\in {\bf B}_R$,
set
$$U_{q,W}=\Phi^{-1}h_{q,W},$$
Since $h_{q,W}$ is holomorphic, $U_{q,W}$ is $J$-holomorphic.
Finally let $\Psi_J$ the map that to $(q,W)\in {\bf B}_R$ 
associates
$$
\big( U_{q,W}(0),{\partial U_{q,W}\over \partial x}(0), \ldots
,{\partial^kU_{q,W}\over \partial^kx}(0)~\big)
\in \R^{2n}\times (\R^{2n} )^k.
$$
If $J=J_{st}$, $\Psi_J$ is the identity mapping. 
If $J$ is close to $J_{st}$ in $\CC^{k,\alpha}$ topology,
it is a small continuously differentiable perturbation of
the identity, whose image therefore contains a neighborhood of
0, with a $\CC^1$ inverse. If $(q,W)=\Psi^{-1}(p,V)$, then
$u_{p,V}=U_{q,W}$ is the desired map. \hfill Q.E.D.
\bigskip

Proposition 1.1 can be rephrased in terms of matching jets. In 
Proposition 1.1'  we rephrase only the part of proposition 1.1
with $\CC^{k,\alpha}$ smoothness assumption.
\bigskip
\noindent {\bf Proposition 1.1'.}
{\it Let $k\ge 1$ and $0<\alpha <1$.
Let $J$ be a $\CC^{k,\alpha}$ almost complex
structure defined in a neighborhood of $0$ in $\R^{2n}$.  If
$\varphi :\D\rightarrow (\R^{2n}, J)$ is a smooth map such that
$\overline\partial_J\varphi (z)= o
(|z|^{k-1})$, then there exists a $J$-holomorphic map $u$ from a \
neighborhood of
$0$ in $\C$ into $\R^{2n}$ such that $|(u-\varphi )(z)|=o
(|z|^{k})$.

If we have a family of maps
$\varphi_t:\D\to (\R^{2n}, J)$ as above, with $\CC^1$
dependence
on $t$ in some neighborhood of $0$ in $\R^\ell$, then
there are $J$-holomorphic maps $u_t$ defined for $t$ near $0$ on a same
neighborhood of $0$ in $\C$ and with $\CC^1$ dependence on $t$,
with $|(u_t-\varphi )(z)|=o
(|z|^{k})$ (uniformly in $t$).}
\bigskip
\noindent {\bf Proof.} Proposition 1.1' follows immediately from
Proposition 1.1, due to the following observation.
Assume that $u$ and $v$ are $\CC^k$ maps from $\D$ into
$(\R^{2n},J)$, where $J$ is an almost complex structure of class
$\CC^{k-1}$. If
$${\partial^\ell u\over \partial^\ell x} (0) =
{\partial^\ell v\over \partial^\ell x} (0),~{\rm for~all}~0\leq l \leq k~,$$
$${\rm and}~\overline\partial_Jv=o(|z|^{k-1})~,~
\overline\partial_Ju=o(|z|^{k-1})~,$$
then $u$ and $v$ have same $k$-jet at 0.
\bigskip
$\overline\partial_Ju=o(|z|^{k-1})$ implies that 
$${\partial u \over \partial y}(z)~=~-J(u(z)){\partial u\over\partial x}(z)~+~
o(|z|^{k-1}).\eqno (1.9)$$
Similarly for $v$. Assume that for some $\ell <k$, we have already shown that
at 0 all the derivatives of $u$ and $v$ of
order $\leq \ell$ coincide. We have to show that for derivatives of
order $\ell +1$, $D={\partial^{\ell +1}\over \partial x^m\partial y^{\ell+1-m}}$,
we have $Du(0)=Dv(0).$

The case $m=l+1$ is covered by the hypothesis. So, take $m\leq \ell~(<k)$, and  write 
$D={\partial\over \partial y}{\partial^\ell \over \partial x^m\partial y^{\ell-m}}$.
Differentiation of $(1.9)$ gives
$$Du=-{\partial^\ell \over \partial x^m\partial y^{\ell-m}}~J(u){\partial u\over \partial x}
+o(1),$$
and similarly for $v$. By the induction hypothesis, and the matching of the pure
$x$-derivatives, one gets $Du(0)=Dv(0)$.\hfill Q.E.D.
\bigskip

\noindent {\bf 1.e.4.} Families of discs.
\smallskip\noindent
Without looking for more generality, we state the next Proposition just as we will 
need it it for proving part $2B$ in Theorem 2.

\noindent {\bf Proposition 1.2.} {\it 
Let $\Omega$ be an open set in $\R^{2n}\simeq \C^n$.
Let $\delta >0$ and $\rho >1$.
Let $\phi_t$ be a family of ($J_{st}$) holomorphic 
discs $\phi_t:\D_\rho \to \Omega$, defined for $-\delta \leq t \leq \delta$,
depending continuoulsy on $t$. Let $\eta >0$. For any almost complex structure
on $\Omega$ close enough to $J_{st}$ in $\CC^{1,\alpha}$ topology on $\Omega$,
there exists a family $\psi_t$ of holomorphic discs
$\psi_t:\D\to \Omega$, continuous in $t$,
 such that for any $t\in [-\delta~,\delta]$ and any $\zeta\in \D$,
$|\psi_t(\zeta)-\phi_t(\zeta)|\leq \eta$.}
\bigskip
\noindent {\bf Proof.} The operators $\Phi_J$ were introduced in and discussed in
the two previous sections. $\Phi_{J_{st}}$ is the identity mapping. Given any compact set
$F$ in $\CC^{1,\alpha}(\overline \D)$, if $J$ is close enough 
$J_{st}$,  $\Phi_J^{-1}$ can be defined on $F$ and is close to the identity.
It maps $J_{st}$ holomorphic discs to $J$ holomorphic discs.\hfill Q.E.D.
\bigskip

\noindent {\bf \S2. Completeness of strictly pseuconvex domains.}
\medskip

\noindent {\bf 2.a. Localization results.} 
\medskip
As a preliminary, we shall start with a first localization (Lemma 2.1) 
that is not needed when dealing with local problems only, for which Lemma
2.2 can be obtained directly. We shall use localization
techniques that are standard in complex manifolds. See in particular 
Proposition 2.1 in [Be], but we avoid any explicit use of
the Sibony metric.
\medskip

\noindent {\bf Lemma 2.1.} {\it Let $D$ be a domain in an almost complex
manifold $(X,J)$, $J$ of class $\CC^1$.  Let $p\in\partial D$.  
Assume that there exists a
neighborhood $U$ of $p$ in $X$ and a continuous function $\rho$ on
$U\cap\overline D$ such that
$$
\left\{\matrix{
\rho (p)=0,\;\;\rho <0\;\hbox{ on $U\cap\overline D-\{p\}$}\hfill\cr
\rho\;\;\hbox{is plurisubharmonic on $U\cap D$.\hskip2.0in}\hfill\cr}
\right.
$$
Then for every $r\in [0,1)$ and for every neighborhood $V$ of $p$ in $X$,
there exists a neighborhood $W$ of $p$ such that if $u:\D\to D$ is a
$J$-holomorphic disc and $u(0)\in W$, then $u(z)\in V$ for
every $z \in\C$ such that $|z|<r$.}

Note that there is absolutely no global assumption made on $D$.
\medskip

\noindent {\bf Proof}.  We first make the function $\rho$ globally defined on $\overline D$ and
identical to $-1$ off a relatively compact subset of $U$ by replacing
$\rho$ by $\max (\kappa\rho ,-1)$ for $\kappa >0$ large enough.  

According to Lemma 1.4, there exists a
plurisubharmonic function $\lambda$ defined near $p$, in $X$, such hat
$\lambda$ is finite and continuous except at $p$ and $\lambda
(p)=-\infty$. 
\medskip

We then replace $\lambda$ by a plurisubharmonic function $\lambda^{\#}$
continuous on $\overline D-\{ p\}$, bounded off any neighborhood of $p$,
plurisubharmonic on $D$ and such that $\lambda^{\#}(q)\to -\infty$ as
$z\to p$.  

Such a function $\lambda^{\#}$ is obtained as follows.
\smallskip

\item{} If $\rho (q)>-{1\over 3}$ set $\lambda^{\#}(q)=C\big(\rho
(q)+{1\over 2}\big)+\lambda (q)$.
\smallskip

\item{} If $-{2\over 3}\le \rho (q)\le -{1\over 3}$ set
$\lambda^{\#}(q)=\max \Big(\rho (q), C\big(\rho (q)+{1\over
2}\big)+\lambda
(q)\Big)$.
\smallskip

\item{} If $\rho (q)\le -{2\over 3}$ set $\lambda^{\#}(q)=\rho (q)$.
\smallskip

If $C$ is taken large enough ($C>0$), the above defines (on $D$)
a global
plurisubharmonic function $\lambda^{\#}$ as desired.

Fix $r\in [0,1)$ and the neighborhood $U$.  Let $u_j:\D\to D$ be a sequence
of $J$-holomorphic discs and
$(\zeta_j)$ a sequence in $\C$, with $|\zeta_j|<r$.  Assume that
$u_j(0)\to p$. We have to show that $u_j(\zeta_j)\to p$.

Since $u_j(0)\to p$, $\rho\circ u_j(0)\to 0$.  By the mean
value property
$$
\int_0^{2\pi}|\rho\circ u_j(e^{i\theta})|{d\theta\over 2\pi}\longrightarrow
0.\eqno(2.1)
$$
For every $N>0$ there exists $\varepsilon >0$ such that
$\rho\circ u_j(z)>-\varepsilon$ implies $\lambda^{\#}\circ u_j<-N$.  Taking into
account that $\lambda^{\#}$ is bounded away from $p$ it follows from
(2.1) that
$$
\int_0^{2\pi} \lambda^{\#}\circ u_j(e^{i\theta}) {d\theta\over
2\pi}\longrightarrow -\infty .
$$
From which it follows that $\lambda^{\#}\circ u_j(\zeta_j)\to -\infty$ and so
$u_j(\zeta_j)\to p$.
\medskip
\hfill Q.E.D. \break
This first localization is followed now by a more precise localization
that will use Lemma 2.1 as a first step.
\medskip

\noindent {\bf Definition}:  We say that the boundary of $D$ (as above) is
{\it strictly\/} $J$-{\it pseudoconvex\/} at $p$.  If there exists a
$\CC^2$ strictly $J$-plurisubharmonic function $\rho $ defined near $p$
in $X$ such that $\nabla\rho \not= 0$ on $\partial D$ (near $p$), and
near $p$, $D$ is defined by $\rho <0$.
\medskip

We equip $X$ with an arbitrary smooth Riemannian metric. 
We still assume $J$ to be of class $\CC^1$.
\medskip

\noindent {\bf Lemma 2.2.} {\it If the boundary of $D$ is strictly
$J$-pseudoconvex at $p_0$, for any $r\in [0,1)$ there exists $\delta >0$
and $C>0$ such that for every $J$-holomorphic disc $u :\D\to D$ with
$\hbox{dist}(u (0), p_0)<\delta$ then
$$
{\rm dist}(u (0), u (z))\le C\sqrt{{\rm dist}(u (0),\partial
D)}\eqno(2.2)
$$
if $|z|<r$.}
\medskip

\noindent
{\bf Proof.} The proof is basically a repetition of the proof of Lemma 2.1 
with
the function $\lambda$ replaced below by the functions $q\mapsto |q-p|^2$.
\smallskip
Fix $r<r_1<1$. Fix a neighborhood $U$ of $p_0$ on which $D$ is
defined by a strictly plurisubharmonic function $\rho$, and diffeomorphic 
to an open set in $\R^{2n}$. On $U$ we will also consider
the Riemannian metric obtained by the identification of $U_0$
with this open set in $R^{2n}$. The distance from $p$ to
$q$ will then be denoted somewhat abusively by
$|p-q|$. 

There is a neighborhood $V\subset U$ of $p_0$, and 
$\varepsilon >0$ (small enough) such that for any $p\in V$,
$q\mapsto \rho_p(q)=\rho (q)-\varepsilon |q-p|^2$ and the functions 
$q\mapsto |q-p|^2$ are $J$-plurisubharmonic on $V$.
Also, for appropriate constants $A$ and $B>0$
$$
-B|q-p|\le \rho_p(q)\le -A|q-p|^2.
$$

By Lemma 2.1 (applied with $\rho_{p_0}$ instead of $\rho$), 
if $u :\D\rightarrow D$ be a $J$-holomorphic map,
and if
$u(0)$ is close enough to $p_0$, $u(z)\in V$ whenever $|z|\leq r_1$. 
Take $p\in \partial D$ such that ${\rm dist}~\big( u(0),\partial D\big)
={\rm dist}~\big( u(0),p\big)$ (so $p\in V$).
For $|z_1|\leq r$, by subharmonicity of $|u(z)-p|^2$ there
is a constant $C$ such that 
$$
|u (z)-p|^2\le C\int^{2\pi}_0 |u (r_2e^{i\theta})-p|^2 {d\theta\over
2\pi}.
$$

\noindent By the mean value property: 
$$
-\rho_p\circ u (0)\ge \int^{2\pi}_0 -\rho_p\circ u (r_2e^{i\theta}) {d\theta\over
2\pi}\ge A\int^{2\pi}_0 |u (r_2e^{i\theta})-p|^2 {d\theta\over
2\pi}\ge {A\over C}|u(z)-p|^2.
$$

It gives 
$$
\big({\rm dist}~(u(z),u(0)\big)^2
\leq \big(\hbox{dist}(u(0), p)+\hbox{dist}(p, u(z)\big)^2\le {\rm C}~\hbox{dist}(u(0),p)
$$
for some other constant $C$, as desired. 
\bigskip
\hfill{Q.E.D.}

\noindent {\bf 2.b. A Schwarz-type Lemma-I and local completeness.}
\smallskip
Lemma 2.1 and 2.2 make possible to work locally if we are interested in the 
distance of a point $p$ of strict pseudoconvexity in the boundary of a domain
$D$, to points in $D$. Take a $\CC^2$-chart 
in a neighborhood $U$ of $p$ such that
$p=0$, $U\cap \overline D\setminus \{ 0\}\subset\{ Z=(z_1,...,z_n)\in \C^n: ||Z||<1, 
{\rm Re}~z_1<0\}$. Using dilations we can assume that $U$ contains the closed unit ball
in $\R^{2n}$ and 
that $J$ is as close as desired to $J_{st}$ in $\CC^1$- 
norm. For any $J$-disc $u:\D\to D$ with $u(0)$ close to $p$, $u(z)\in U$ 
for $|z|<{1\over 2}$. So due to (1.5) and (2.2) we see that 

$$
|\nabla u(z)|\le C\sqrt{\hbox{dist}(u(0),\partial 
D)}\le C\sqrt{-Reu_1(0)}\quad\hbox{ if $|z|\le 1/4$}. 
$$
We state the next Lemma so that, starting form above, it can be used
without (the obvious) rescaling (done in the proof).

\noindent {\bf Lemma 2.3 (Schwarz-type Lemma-I).}  {\it Let $J$ be of class 
$\CC^1$ on the closure of the unit ball in $\R^{2n}$.
Assume that there exists $A>0$ such that for every 
$J$-holomorphic disc $u:\D_{1\over 2}\to U^-:=\{Z=(z_1,...,z_n)\in \C^n:~|Z|<1,~
{\rm Re}~z_1<0\}$ such that 
$u(0)$ is close enough to $0$:
$$
|\nabla u(z)|\le A\sqrt{-{\rm Re}~u_1(0)}\eqno(2.3)
$$
if $|z|\le 1/4$.
Then there exists $C=C(A)>0$ such that for every $J$-holomorphic $u:\D\to U^-$ 
with $u(0)$ close enough to $0$ one has
$$
|\nabla {\rm Re}~u_1(0)|\le C|{\rm Re}~u_1(0)|.\eqno(2.4)
$$
}
\medskip

\noindent {\bf Proof}. Of course, replacing the function $u(z)$ by the function
$u({z\over 4})$, we can assume that $u$ is instead defined on $\D$ and
that $(2.3)$ holds for $|z|<1$.

Formula (1.2) applied to the coordinate function 
$\lambda = {\rm Re}~z_1$ gives 
$$
\Delta {\rm Re}~u_1(z)=(dd^c_{J}z_1)_{y(z)}\Big({\partial u\over\partial 
x}(z), J{\partial u\over\partial x}(z)\Big).\eqno(2.5)
$$
Using the estimate (2.3) one gets
$$
|\Delta {\rm Re}~u_1(z)|\le C_1|{\rm Re}~u_1(0)|\eqno(2.6)
$$
if $|z|\le 1$, for a conctant $C_1$ depending only on $J$ and $A$.
The estimate (2.4) follows from $(2.6)$ and from the condition
${\rm Re}~u_1<0$.
\bigskip
So, set $f(z)={\rm Re}~u_1(z)$.

By $\Delta f~\tilde{}$, we denote the function
equal to $\Delta f$ on $\D$, and 0 elsewhere. With $C_1$ as above,
set $$g(z)=f(z)-[{1\over 2\pi}\ln |\zeta|*\Delta f~\tilde{}~](z)-C_1|f(0)|.
\eqno (2.7)$$
From $(2.6)$, we get $|{1\over 2\pi}\ln |\zeta|*\Delta f~\tilde{}~|\leq C_1|f(0)|$.
So $g<0$ and $|g(0)|< (2C_1+1)|f(0)|$.

The classical Schwarz Lemma for negative harmonic functions gives:
$|\nabla g(0)|\leq 2 |g(0)|$. Therefore (2.7) yields;
$$|\nabla f(0)|\leq |\nabla g(0)|+{\rm Sup}~|\Delta f~\tilde{}~|
\leq 2|g(0)|+C_1|f(0)|\leq (5C_1+2) |f(0)|.$$
\hfill Q.E.D.
\bigskip
The localization Lemma 2.1 and 2.2, and  the above allow us to apply
Lemma 1.1 to the function $\chi={\rm Re}~z_1$ in order to immediately get:
\medskip
\noindent {\bf Proposition 2.1.}  {\it Let $D$ be a domain in an almost 
complex manifold $(X,J)$, $J$ of class $\CC^1$.  Let $p\in\partial D$.  If 
the boundary of $\partial D$ is strictly $J$-pseudoconvex at $p$, the 
point $p$ is at infinite Kobayashi distance from the points in $D$.}
\noindent

Note that there is still no global hypothesis on $D$. 
So $D$ may very well not be hyperbolic i.e. the Kobayashi-Royden
`distance' may be only a pseudo-distance in the interior. An easy example
are obtained by the blow up of a point. It has been pointed out  
at the end of 1.c that in an almost complex manifold, every point has a basis of
hyperbolic neighborhoods. Lemmas 1.3 and 1.5 show that the small balls
(with respect to any given Riemannian metric) are complete hyperbolic,
and so are the intersections of small balls with strictly pseudoconvex
domains. 
\bigskip\noindent
{\bf 2.c. End of the Proof of Theorem 1.}
\bigskip To get Theorem 1, hyperbolicity of the domain, in case $D$
does not contain complex lines, is the only point left. 

It follows from the Brody reparametrization lemma. 
If there
is a sequence of points $p_j\in D$ and  unit (with respect to some metric) 
vector $v_j\in T_pX$ such that $|v_j|_k$ tends to 0, then
there is a sequence of $J$-holomorphic maps $u_j:\D_{R_j}\to D$ such 
that $du_j({\partial\over \partial x})=v_j$ and $R_j\to \infty $. After
reparametrization we get a sequence $\tilde u_j:\D_{R_j}\to D$, still
$J$-holomorphic and such that $\sup \{ |d\tilde u_j(z)\cdot {R_j^2-|z|^2
\over R_j^2}|:z\in \D_{R_j}\} =
|d\tilde u_j(0)|=1$. A subsequence converges to a non-constant map from $\C$
into $\overline D$, which due to the strict pseudoconvexity of the boundary of
$D$ must be a map from $\C$ into $\D$. 

\bigskip Theorem 1 is proved.

\bigskip

\noindent {\bf \S3. Estimate of a Calderon-Zygmund Integral and Schwarz-type
lemma.}
\bigskip

Fix the following function $\phi (r)=r\ln {1\over r}, r>0$.
We have already introduced the  class $\CC^{1,\phi}$.
Let us 
introduce more generally sub-Lipschitzian
classes $\CC^{k,\phi}$ which are the spaces of functions or maps
$ f\in \CC^k$ that locally satisfy  
${|f^{(k)}(z')-f^{(k)}(z)| \over \phi (|z'-z|)}\leq C<\infty$ (for some
positive $C$).
We also define the Banach space $\CC^{k,\phi}(\overline \D)$ as the
space of complex valued functions on $\D$ that staisfy:
$$
\|f\|_{\CC^k(\overline \D)}+\sup_{|z-z|\leq {1\over 2}} {|f^{(k)}(z')-f^{(k)}(z)|
\over 
\phi (|z'-z|)}<\infty ,$$
with the left hand side defining the norm $\|\cdot \|_{k,\phi}$, $k\ge 0$. 
$\CC^{0,\phi }$ will be denoted simply by $\CC^{\phi }$.

\noindent {\bf Lemma 3.1.} {\it There exists a constant $C$ such that for
all complex-valued functions $f\in \CC^{\phi}(\R^2)$, and $g\in \CC^{1,\phi}(\R^2)$ such that:

\item{(1)} $\| f\|_{\CC^{\phi}(\overline\D )}\leq 1$,  $\|
g\|_{C^{1,\phi}(\overline\D )}\le 1$;

\item{(2)} $\| g\|_{C^0(\overline\D )}\le {1\over 2}$;

\item{(3)} $|f(z)|\le |g(z)|$ and $g(z)\not= 0$ for all $z\in\D$,

\smallskip\noindent 
one has
$$
\left|\int_{\overline\D} {1\over z^2} {f(z)\over g(z)} dzd\overline
z\right| \le C\cdot\log {1\over |g(0)|}.\eqno(3.1)
$$
\smallskip

}
\smallskip
\noindent {\bf Proof}.  We can assume $g(0)>0$. Set $\delta =g(0)$ and split 
$$
\int_{|z|\le 1} {1\over z^2} {f(z)\over g(z)} dxdy=\int_{|z|\le\delta /4}
{1\over z^2} {f(z)\over g(z)} dxdy+\int_{\delta /4\le |z|\le 1} {1\over
z^2} {f(z)\over g(z)} dxdy,
$$
where clearly
$$
\int_{\delta /4\le |z|\le 1}\left| {1\over z^2} {f(z)\over g(z)}\right|
dxdy\le C\cdot \log {1\over\delta},
$$
since $\delta =g(0)\le {1\over 2}$.
In order to estimate $|\int_{|z|\le\delta /4}
{1\over z^2} {f(z)\over g(z)} dxdy|$,  we shall use the following
cancellations (with the first integral in the sense of principal value):
$$
\int_{|z|\le\delta /4} {1\over z^2}dxdy =\int_{|z|\le\delta /4} {z\over z^2}dxdy
=\int_{|z|\le\delta /4}{\overline z\over z^2}dxdy =0.
\eqno(3.2)
$$
Write $f(z)=f(0)+ R_1(z)$, and $g(z)=g(0)+Az+B\overline
z+Q_2(z)$.  Due to the condition of the Lemma that $\| f\|_{\phi}, \|
g\|_{1,\phi}\le 1$ we have (with appropriate definitions of the norms)
$$|A|, |B|\leq 1,~ |R_1(z)|\leq |z|\ln {1\over |z|},~|Q_2(z)|\le |z|^2
\ln {1\over |z|}. 
$$ 
For $|z|\leq {\delta\over 4}$, we have that $g(z)=\delta
(1+{Az+B\overline z\over\delta}+{Q_2(z)\over\delta})$ with
$|{Az+B\overline
z\over\delta}|+|{Q_2(z)\over\delta}|\le {1\over 4}+{1\over
4}+{\delta\ln {1\over \delta }\over 16}\le {2\over 3}$.  

For dealing with $1\over g$, we shall simply use that if $|a|+|b|\leq {2\over 3}$,
${1\over 1+a+b} = 1-a+r$ with $r\leq C (|a|^2+|b|)$, for some universal
constant C. We will apply it with $a= {Az+B\overline z\over\delta}$,
and $b={Q_2(z)\over\delta}$.

So, we can write that
$$
{f\over g} ={1\over\delta }[f(0)+R_1(z)]\left[
1-\left({Az+B\overline z\over\delta}\right) +S_2(z)\right],
\eqno(3.3)
$$
with $|S_2(z)|\le C\left(\left|{Q_2(z)\over \delta }\right|+\left|{Az+B\bar z\over
\delta }\right|^2\right)\le C({|z|^2\over \delta}\ln {1\over |z|}+{|z|^2\over
\delta^2})$.  Write (3.3) as
$$
{f\over g}={f(0)\over\delta}-f(0){A z+B\overline z\over
\delta^2}+T(z),
$$
where $|T(z)|={f(0)+R_1(z)\over\delta} S_2(z)+ {R_1(z)\over \delta }
(1-{Az+B\overline z\over \delta})$.

\medskip
\item{(1)} From $|{f(0)|\leq |g(0|=\delta}$,
and $|R_1(z)|\leq |z|\ln {1\over |z|}$, we
get that $\left|{f(0)+R_1(z)\over\delta} S_2(z)\right|$ $\le C(1+|z|\ln {1\over 
|z|})({|z|^2\over\delta^2}+{|z|^2\over \delta}\ln{1\over |z|})\le 
C({|z|^2\over \delta^2}+{|z|^2\over \delta}\ln{1\over |z|}+{|z|^3\over \delta}
\ln^2{1\over |z|})$.

\item{(2)} $\left| R_1(z) (1-{Az+B\overline z\over \delta}) \right| \leq {3|z|}
\ln {1\over |z|}$,
~so ${1\over \delta} \left| R_1(z) (1-{Az+B\overline z\over \delta})\right| 
\leq {3|z|\over \delta}\ln{1\over |z|}.$

\smallskip
Using cancellation properties of the Calderon-Zygmund operator and the
claim we have proved one gets
$$
\left|
\int_{|z|\le{\delta\over 4}} {1\over z^2} {f(z)\over g(z)} dxdy\right|
=\left| \int_{|z|\le {\delta\over 4}} {T(z)\over z^2} dxdy\right|\le
C \int_{|z|\le\delta /4} {1 \over \delta |z|}\ln{1\over |z|} dxdy\le
C\ln{1\over \delta}.
$$
\hfill{Q.E.D.}

\medskip
Lemma 3.2 (in which we focus on the behavior near the puncture)
will be a generalization of the
standard Schwarz Lemma which gives the following estimate for a holomorphic
map $g$ from the unit disc into the punctured unit disc:
$$
|g'(0)|\le 2|g(0)| \log \left|{1\over g(0)}\right|.
$$

\smallskip

\noindent {\bf Lemma 3.2 (Schwarz-type Lemma-II).} {\it For $A$ and $B>0$ there exists $C>0$ such
that for every map $g$ from the unit disc into the punctured disc
$\D_{1/2}-\{ 0\}$ satisfying
$$
\| g\|_{1,\phi}\le A\;\;\hbox{ and }\;\; \left|{\partial g(z)\over\partial
\overline z}\right|\le B|g(z)|
$$
for all $z\in D$, one has
$$
|\nabla g(0)| \le C|g(0)|\log {1\over |g(0)|}.
$$
\smallskip

}

\noindent {\bf Proof}.  Extend ${-{\partial g\over\partial \overline z}\over g}$ 
to $\C$ (identified with $\R^2$), by
setting it to be equal to $0$ outside the unit disc.  Set
$$
w={1\over \pi z}\; *\; {-{\partial g\over\partial \overline z}\over g}.
$$
\smallskip

We have $|w|\le 2 B$, and on the unit disc,
$$
{\partial w\over\partial \overline z}={-{\partial g\over\partial \overline z}
\over 
g}.
$$
So $h=ge^w$ is holomorphic.  The holomorphic function $h$ never vanishes,
and it takes values in the disc of radius ${1\over 2}e^{2B}$.  So, the Schwarz
Lemma for holomorphic function gives us a bound:
$$
|\nabla h(0)|\le 2|h(0)|\log \left| {C_1\over h(0)}\right|,
$$
for some constant $C_1$ depending only on $B$. We have
$$
\nabla g (0)=e^{-w(0)} \nabla h(0)+h(0)e^{-w(0)}\nabla w(0).
$$
Note that we have a bound for $|w(0)|$, so $h(0)$ and $g(0)$ are
comparable (bounded ratios), and we can use the above estimate for $\nabla
h(0)$.  All what is left is to have a correct estimate of $\nabla w(0)$
(by a multiple of $\log {1\over |g(0)|}$).  The $\overline z$ derivative
of $w$ is simply ${-{\partial g\over \partial \overline z}\over g}$, which has
modulus $\le B$.  The needed estimate of ${\partial w\over\partial z}(0)$
which is given by the integral (defined as a principal value):
$$
\int_{\D} {1\over z^2} {-{\partial g\over\partial \overline z}\over g}~,
$$
is given by Lemma 3.1, applied to the function 
 $g_1: z\mapsto g({1\over k} z)$ instead of $g$, with
$k\leq {\rm max}~ (1, \|g\|_{\CC^{1,\phi}})$, and $f=c
{\partial g_1\over\partial\overline z}$, for an appropriate constant
$c$.
\bigskip
\hfill{Q.E.D.}
\bigskip

\noindent {\bf \S4. Distance to a complex hypersurface.}
\smallskip\noindent\sl Proof of (2.A). \rm 
Due to hyperbolicity the question can be
localized.  Therefore we will suppose that $X$ is a neighborhood of $0$ in
$\R^{2n}$ and that $M=\{ Z=(z_1,\ldots , z_n)\in X: z_n=0\}$, where
$z_j=x_j+iy_j$ are complex coordinates in $\R^{2n}\cong \C^n$, unrelated
to the almost complex structure $J$. Since
$M$ is of class $\CC^{3}$, in the new coordinates we can keep
$J$ of class $\CC^{2}$. Without loss of generality we
can assume that $J(0)=J_{st}$, the standard complex structure on $\C^n$.
We shall write $Z=(z_1,\ldots , z_n)=(Z'; z_n)$ with $Z'=(z_1,\ldots ,
z_{n-1})$.  Because $z_n=0$ is a $J$-complex hypersurface we already have
the almost complex structure $J$ given along $z_n=0$ by
$$
J(Z',0)=\left(\matrix{A(Z')&\alpha\cr
0 &\beta\cr}\right),
$$
where $\alpha$ is a $(2n-2)\times 2$ matrix while $\beta$ is a $2\times 2$
matrix.  After shrinking the neighborhood of $0$ if needed, consider the
$\CC^{2,\alpha}$  change of variables given by
$$
\Phi (z_1,\ldots , z_n)=(Z';0)+(0; x_n,0)+y_nJ(Z';0){\partial\over\partial
x_n},
$$
(using obvious identification of $\R^{2n}$ and $T\R^{2n}$).
In the new coordinate system we will have
$$
J(Z';0)=\left(\matrix{
A(Z')&0\cr
0&J^{(2)}_{st}\cr}\right)\eqno(4.1)
$$
where $J_{st}^{(2)}=\left(\matrix{0&-1\cr 1&0\cr}\right)$.

From now on we will work with coordinates in which (4.1) holds,
with $J$ now of class $\CC^{1}$.

Let $u :\D\to (\R^{2n}, J)$ be a $J$-holomorphic map. Recall the 
condition for
$J$-holomorphicity in the form (1.7)

$$
{\partial u\over\partial\overline z}+Q_J(u) {\partial u\over\partial
z}=0\eqno(4.2)
$$
with $Q_J(u)=[J(u)+J_{st}]^{-1}[J(u)-J_{st}]$, where
${\partial\over\partial \overline z}$ and ${\partial\over\partial z}$
refer to the standard complex structures on $\D$ and $\C^n$, i.e.\ in
$\R^{2n}$ coordinates
$$
{\partial u\over\partial\overline z}={1\over 2}\left({\partial
u\over\partial x}+J_{st}{\partial u\over\partial y}\right), \;\;
{\partial u\over \partial z}={1\over 2}\left( {\partial u\over\partial
x}-J_{st}{\partial u\over\partial y}\right).
$$
Due to (4.1) 
$$Q_J(Z',0)=\left(\matrix{B(Z')&0\cr
0&0\cr}\right).
$$
The Cauchy-Riemann equation (4.2) for $J$-holomorphic maps $u:\D\to
(\R^{2n},J)$ gives for the component $u_n$:
$$
{\partial u_n\over\partial\overline z}=\sum^n_{k=1} b_k(u_1,\ldots
, u_n) {\partial u_k\over\partial z}
$$
with $b(u_1,\ldots ,u_{n-1},0)\equiv 0$.  For some constant $C$ we
therefore get an inequality
$$
\left|{\partial u_n\over\partial\overline z}\right|\le C|u_n| \sum^n_{k=1}
\left| {\partial u_k\over\partial z}\right|.\eqno(4.3)
$$
On the right-hand side $\left|{\partial u_k\over\partial z}\right|$ will
be bounded by the $C^{1,\phi}$ regularity of $J$-holomorphic maps, 
the factor $|u_n|$ is the important feature.

Since the problem is purely local and since
(after the changes of variables), the almost complex
structure $J$ is of class $\CC^{1}$,  we can use
shrinking and rescaling in order to assume that $(X,J)=(\D^n, J)$ with
$\D^n$ the unit polydisc in $\C^n\cong \R^{2n}$, equipped with an
almost complex structure $J$ sufficiently close to the standard structure
in $\CC^1$ topology. Then, by Remark 3 in 1.d,
$J$-holomorphic maps from $\D$
into $(\D^n, J)$ have uniformly bounded $C^{1,\phi}$ norm on $\D_{1/2}$.

From some (other) constant $C$ we
therefore get:
$$
\left|{\partial u_n\over\partial \overline z}\right|\le C|
u_n|\;\;\hbox{ on $\D_{1/2}$},
$$
and
$$
\| u_n\|_{\CC^{1,\phi}(\D_{1/2})}\le C.
$$
Only now we are going to use the hypothesis that the map $u$ should avoid
the hyperplane $M=\{ z_n=0\}$, i.e.\ $u_n\not= 0$.  We apply Lemma 3.2 to
the restriction of ${u_n\over 2}$ to $\D_{1/2}$ (i.e.\ to
$g(z)={1\over 2} u_n\left({z\over 2}\right)$), to get for some (other)
constant $C$:
$$
|\nabla u_n(0)|\le C|u(_n0)|\log {1\over 2|u_n(0)|}.
$$
By Lemma 1.1, if $p=(p_1,\ldots , p_n)\in \D^n$ and $p_n\not= 0$,
the Kobayashi distance $d_{\D^n\setminus M} (p,M)$ from $p$ to $M$,
relative to $\D^n\setminus M$ and $J$ is infinite.
\smallskip\noindent
This achieves the proof of (2.A).

\noindent {\bf Proof of $(2.B)$.} Since $M$ is not $J$ complex at $p$ and is of class $\CC^3$, 
we find a $\CC^3$ change of coordinates so that: $p=0$, $J(0)=J_{st}$,  
$M$ coincides with $\R^2\times\C^{n-2}$ on 
the  ball of radius 6, and that this ball is included in $D_1$. Note that $J$ stays of class $\CC^2$.
Consider the following family of $J_{st}$ holomorphic discs defined for
$|\zeta |\leq 1$, $|t|\leq 1$:
$$\phi_t(\zeta)=(\zeta , t+i\zeta^2,0,\ldots ,0),~{\rm if}~t\geq 0,$$
$$\phi_t(\zeta)=(\zeta , i{t\over 8}\zeta +i\zeta^2,0,\ldots ,0)~{\rm if}~t< 0.$$
This family has the following properties.
\smallskip\noindent
(a) $\phi_t(\D\setminus \D_{1\over 2})\subset K$ for some compact $K$ not intersecting $M$.

\noindent
(b) $\phi_t(\bar \D)\cap M=\emptyset $ for $t>0$;

\noindent
(c) $\phi_0(\bar \D)\cap M=\{ 0\} $;

\noindent
(d) $\phi_t(\bar \D)$ intersects $M$ transversally for $t<0$.
\bigskip
Using dilations, replacing $J(z)$ by $J(\epsilon z)$, we can assume that
$J$ is as close as we need to $J_{st}$. Note that these dilations
leave $M=\R^2\times\C^{n-2}$ invariant, but that we do not rescale
the family $(\phi_t)$. Applying Proposition 1.2 we then get a family
of $J$ holomorphic discs $\psi_t$ with the following properties:
\smallskip\noindent
(a') $\psi_t(\{|z|={3\over 4}\} )\subset K$ for some compact $K$ not intersecting $M$.

\noindent
(b') $\psi_t(\bar \D)\cap M=\emptyset $ for $t>\delta_0$, for some $\delta_0$ close to
0.

\noindent
(c') $\psi_{\delta_0} (\bar \D)\cap M\neq \emptyset $.
\bigskip
Clearly any point in  $p'=\psi_{\delta_0} (\zeta_0)\in M$ is at finite distance
from points in the complement of $M$. Indeed, for $t>\delta_0$ the distance from
$\psi_t({3\over 4})$ to $\psi_t(\zeta_0)$ in $D\setminus M$ stays bounded as 
$t\to \delta_0$, but $\psi_t(\zeta_0)$ tends to a point in the complement of $M$.

\hfill Q.E.D.

\bigskip\noindent {\bf \S5. Hyperbolic distance to a real hypersurface. }
\smallskip

\noindent\bf
5.a Construction of the family of discs.

\smallskip\rm

From Proposition 1.1' we will deduce:
\medskip

\noindent {\bf Proposition 5.1.} {\it Let $J$ be a $\CC^{2,\alpha}$
almost complex
structure defined near $0$ in $\R^{2n}$.  Let $M$ be a $\CC^2$
hypersurface in $(\R^{2n}, J)$ defined by $\rho =0$, $0\in M$,
$\nabla\rho (0)\not= 0$.  If there is a complex tangent vector $Y\in
TM(0)$ such that $dd_J^c\rho (Y,J(0)Y)> 0$ then there is a family of
$J$-holomorphic embedded discs $u_t :z\to\varphi_t(z)\in\R^{2n}$
($z\in\D$) that depend continuously on $t\in [0,1]$ such that
$$\left\{
\matrix{\hbox{if $t>0$, $u_t(\D )\subset \{ \rho >0 \}$ 
}\hfill\cr
0=u_0(0)=u_0(\D )\cap M.\hfill\cr}\right.
$$
\smallskip

}

\noindent {\bf Proof}. Let $\psi :\D_R\to (\C^{2n},J)$ be some $J$-
holomorphic disk with $\psi (0)=0$ and ${\partial \psi \over \partial x}(0)
=Y$. Then $\psi \in \CC^{3,\alpha }$ and after a $\CC^{3,\alpha}$
change of coordinates we can assume that $\psi
:z\mapsto (\hbox{Re }z, \hbox{Im }z, 0,\ldots , 0)$ 
(($z,0, \ldots ,0)$, when using complex coordinates)
is a $J$-holomorphic
disc, $Y=(1,0,\ldots , 0)$, $J(0)Y=(0,1,0,\ldots , 0)$ and $J=J_{st}$
along $\R^2\times \{ 0\}$. 
Moreover we can choose coordinates so that the
vector $(0,0,1,0,\ldots , 0)~(\in \R^{2n})$ is a normal vector to $M$ at $0$ with
$\nabla\rho (0)[(0,0,1,0,\ldots , 0)]=1$.
In the new coordinates $J$ is still of class $\CC^{2,\alpha}$.
\medskip

Since $\psi$ is a $J$-holomorphic disc $d^c(\rho \circ \psi )=\psi_*d_J^c\rho$.
By differentiation $dd^c(\rho \circ \psi )=\psi_*dd_J^c\rho
\not= 0$.  Therefore the Taylor
expansion of $\rho \circ \psi$ at $0$ is
$$
\rho \circ \psi (z) =\hbox{Re }az^2+b|z|^2+o (|z|^2),
$$
with $b> 0$.  Using complex $\C^n$ coordinates in the right-hand side,
consider $\varphi (z)=(z,-az^2,0,\ldots , 0)$.  It need not be a 
$J$-holomorphic map. But since $J=J_{st}$ along $\C\times \{ 0\}$ and since
$\varphi$ is holomorphic for the standard structure, we have
$\overline\partial_J \varphi =O (|z|^2)$.  By Proposition 1.1 there exists a
germ of $J$-holomorphic disc $u_0$ such that, still using complex
coordinates,
$$
u_0(z)=(z, -az^2,0,\ldots , 0)+o (|z|^2).
$$
It is immediate to check that
$$
\rho \circ u_0(z)=b|z|^2+o (|z|^2).
$$
Then one has just to restrict $u_0$ to a small disc to be identified
with $\D$.  We then have $u_0(0)=0$ but
$b\big( \rho\circ u_0(z)\big)~>~0$ if $z\not= 0$.

The construction of the discs $u_t$ is exactly similar replacing
$M=\{ \rho =0\}$ by positive level sets of $\rho$. 

\bigskip
\hfill{Q.E.D.}

\medskip
\noindent {\bf 5.b. Proof of Theorem 3}.
After change of variables $p=0$, we are in the situation of Proposition
5.1.  Then, the sequence of points $u_{1/n}(0)$ tends to $0$, their
respective Kobayashi distance to $u_{1/n}(1/2)$ in the complement of
$M$ stay bounded, but the points $u_{1/n}(1/2)$ stay in a compact
subset of the complement of $M$. This establishes Theorem 3.
\medskip

\noindent {\bf Proof of a stronger version of (2.B)}. 
We now assume that $J$ is of class $\CC^{2,\alpha}$, and we wish
to show that in $(2.B)$ one can take $p=p'$.
Let $M$ be a real codimension $2$ submanifold of class $\CC^2$ in $(\R^{2n}, J)$, $0\in M$.
Assume that $Y\in TM(0)$ but $J(0)Y\not\in TM(0)$.  Consider a
hypersurface $\widetilde M$ defined by $\rho =0$ containing $M$ such that
$J(0)Y$ is tangent to $\widetilde M$ at $0$.  It is an easy exercise to
show that by adding a quadratic term to $\rho$ if needed (therefore not changing the
tangent space) but keeping $\rho =0$ on $M$ we can get
$$
dd^c_J\rho (Y,J(0)Y)\not= 0.
$$
Observe that if $g$ and $h$ are functions defined near $0$, with $h(\zeta
)=O (|\zeta |^2)$ and $J(0)=J_{st}$ then
$dd^c_J(g+h)(0)=dd^c_Jg(0)+dd^ch(0)$.  So it is enough to take $h=0$ on
$M$ but $dd^ch(Y,J(0)Y)\not= 0$, ($d^c=d^c_{J_{st}}$).

This being done, $0$ is at finite distance from points in the complement
of $\widetilde M$, for the Kobayashi distance relative to the
complement of $\widetilde M$, and thus a fortiori for the Kobayashi
distance relative to the complement of $M$.
\bigskip
\noindent {\bf Remark.} Any $\CC^2$ submanifold $M$ of codimension $>2$
is a submanifold of a submanifold of codimension 2,
whose tangent space (at any chosen point of $M$) is not $J$-complex. 
So, as above, points of $M$ are at finite distance from points
in the complement of $M$.
But a more direct argument can be given, simply based on the fact that,
by simple count of dimensions, $J$-holomorphic discs generically miss
$M$, and $\CC^2$ regularity of $M$ is not needed.
\bigskip
\noindent {\bf \S6. An example.}
\smallskip

\bigskip\noindent
{\bf 6.a. $d^c_J$ and Levi foliation.}
\bigskip
If $M$ is a real hypersurface defined by $\rho
=0$, with $\nabla\rho\not= 0$, a tangent vector $Y$ to $X$ at a point
$p\in M$ is a complex tangent vector to $M$ if and only if $d\rho
(Y)=d^c_J\rho (Y)=0$.

The question of foliation of $M$ by $J$-complex hypersurface is much the
same as in the complex setting.  The question is to know whether when $Y$
and $T$ are complex tangential vector fields to $M$, the Lie bracket
$[Y,T]$ is also complex tangential i.e.\ if $d_J^c\rho [Y,T]=0$.  By the
definition of $d$, $dd_J^c\rho (Y,T)=Yd_J^c\rho (T)-
Td_J^c\rho (Y)-d_J^c\rho ([Y,T])= -d^c_J\rho ([Y,T])$.  Thus the Frobenius
condition $d_J^c\rho ([Y,T])=0$ is equivalent to $dd^c_J\rho (Y,T)=0$.

\smallskip
\noindent {\bf Remark}.  It is well known that for complex manifolds, if a
hypersurface is not Levi flat, its complement is never complete hyperbolic.

There is a very simple and fundamental difference between the complex case
and the almost complex case.  Consider the two conditions on a
hypersurface
\smallskip

\item{(1)} For every complex tangent vector field $Y$, $[Y,JY]$ is complex
tangent (i.e.\ $dd^c\rho (Y,JY)=0$).

\item{(2)} For every complex tangent vector fields $Y$ and $T$, $[Y,T]$ is
complex tangent.
\smallskip

For almost complex manifolds of real dimension $>4$, (1) does not imply (2)
(example above). But for complex manifolds it does.  The easiest way to
see it is by considering the Levi form as a hermitian form on the complex
tangent space, or (just a different writing) by using
$dd^c=2i\partial\overline \partial$.
Here we sketch a direct argument in terms of the Lie brackets.
If (1) holds $[Y+JT, JY-T]$
is complex tangent.  But $[Y+JT, JY-T]=[Y,JY]+[JT,J(JT)]- [Y,T]+[JT,JY]$.  If
$J=J_{st}$, $[JY, JT]=[Y,T]$ modulo a complex tangent vector field, as follows
from the vanishing of the Nijenhuis tensor (for complex manifolds):
$[Y,T]+J[JY,T]+J[Y,JT]-[JY,JT]=0.$  Hence $[Y,T]$ is complex tangent.

For almost complex manifolds of real dimension $4$, (1) and (2) are
obviously equivalent since if $Y\not= 0$, $Y$ and $JY$ generate the
complex tangent space.

\bigskip

\noindent
{\bf 6.b.} In an almost complex manifold of real dimension $4$, the complement of a
hypersurface is (locally) complete hyperbolic if and only if that
hypersurface is foliated by $J$-complex curves.  We now present the following example on $\R^6$, (in which it is to
be noticed that we do not need to restrict to bounded regions).  We use
coordinates $(x_1, y_1, x_2, y_2,x_3,y_3)$.  On $\R^6$ we define the
vector fields
$$
L_1={\partial\over\partial x_2}+x_1 {\partial\over\partial x_3},\quad
L_2={\partial\over\partial y_2}-y_1 {\partial\over\partial x_3}
$$
and we define the almost complex structure $J$ by setting:
$$\eqalign{
J\left({\partial\over\partial x_1}\right) &={\partial\over\partial
y_1}\qquad\hbox{(so $J\left({\partial\over\partial y_1}\right)
=-{\partial\over\partial x_1}$)}\cr
&\cr
J(L_1)&=L_2\qquad\hbox{(so $J(L_2)=-L_1$)}\cr
&\cr
J\left( {\partial\over\partial x_3}\right) &={\partial\over\partial
y_3}\quad \hbox{(so $J\left({\partial\over\partial
y_3}\right)=-{\partial\over\partial x_3}$).}\cr}
$$
Note that the functions $z_1=x_1+iy_1$ and $z_2=x_2+iy_2$ are 
$J$-holomorphic but that $z_3=x_3+iy_3$ is not.  Finally we simply consider
the hypersurface $M=\{ y_3=0\}$.  At each point the tangent space to $M$
is generated by ${\partial\over\partial x_1}$, ${\partial\over\partial
y_1}$, $L_1$, $L_2$ and ${\partial\over\partial x_3}$.  The complex
tangent space ($TM\cap JTM$) is therefore spanned by
${\partial\over\partial x_1}$, ${\partial\over\partial y_1}$, $L_1$,
$L_2$.  Since $\left[{\partial\over\partial x_1}, L_1\right]
={\partial\over\partial x_3}$ is not a complex tangent vector, $M$ is not
foliated by $J$-complex hypersurfaces.  However:
\medskip

\noindent {\bf Proposition 6.1.} {\it For any point $p\in \R^6\setminus M$
the Kobayashi pseudo-distance from $p$ to $M$ is infinite.}
\medskip

It therefore follows from Theorem 2 that for any complex tangential vector
field $Y$ on $M$, $[Y,JY]\in TM\cap JTM$.  It can easily be checked
directly, but we won't need it.
\medskip

\noindent {\bf Proof}.
We will need
$$\eqalign{
J\left({\partial\over\partial x_2}\right) &=J\left( L_1-x_1
{\partial\over\partial x_3}\right) ={\partial\over\partial
y_2}-y_1{\partial\over\partial x_3}-x_1{\partial\over\partial y_3}\cr
&\cr
J\left({\partial\over\partial y_2}\right) &=J\left( L_2 + y_1
{\partial\over\partial x_3}\right) =-{\partial\over\partial
x_2}-x_1{\partial\over\partial x_3}+y_1{\partial\over\partial y_3}.\cr}
$$
We now write in detail the condition in order that a map $u
:\D\to (\R^6, J)$ be $J$-holomorphic.  In $\D$ we use coordinate
$z=x+iy$, and we write $u =(X_1, Y_1, X_2, Y_2, X_3, Y_3)$.

The condition for $J$-holomorphicity is ${\partial u\over\partial
y}=J{\partial u\over\partial x}$.  It gives:
$$\eqalignno{
&\left.\matrix{{\partial X_1\over\partial y}=-{\partial Y_1\over\partial
x}\cr
\cr
{\partial Y_1\over\partial y}={\partial X_1\over \partial
x}\cr}\right\}&(6.1)\cr
&\cr
&\left.\matrix{{\partial X_2\over\partial y}=-{\partial Y_2\over\partial
x}\cr
\cr
{\partial Y_2\over\partial y}={\partial X_2\over \partial
x}\cr}\right\}&(6.2)\cr
&\cr
&\left.\matrix{{\partial X_3\over\partial y}=-{\partial Y_3\over\partial
x}-Y_1 {\partial X_2\over\partial x}-X_1{\partial Y_2\over \partial x}\cr
\cr
{\partial Y_3\over\partial y}={\partial X_3\over \partial
x}-X_1{\partial X_2\over\partial x}+Y_1{\partial
Y_2\over\partial x}\cr}\right\}&(6.3)\cr}
$$
(6.1) and (6.2) merely say that $Z_1=X_1+iY_1$ and $Z_2=X_2+iY_2$ are
holomorphic functions of $z=x+iy$. 

We now compute $2{\partial\over\partial\overline z}Z_3$:
$$\eqalign{
2{\partial\over\partial\overline z} (Z_3)&=\left({\partial X_3\over
\partial x}-{\partial Y_3\over \partial y}\right) +i\left({\partial
X_3\over\partial y}+{\partial Y_3\over\partial x}\right)\cr
&\cr
&=\left( X_1{\partial X_2\over\partial x}-Y_1{\partial Y_2\over\partial
x}\right) +i\left(-Y_1 {\partial X_2\over\partial x}-X_1{\partial
Y_2\over\partial x}\right)\cr
&\cr
&=(X_1-iY_1) {\partial X_2\over\partial x}-(Y_1+iX_1) {\partial
Y_2\over\partial x}\cr
&\cr
&=\overline Z_1 {\partial X_2\over\partial x}-i\overline Z_1{\partial
Y_2\over\partial x}=\overline Z_1 {\partial \overline Z_2\over\partial
x}.\cr}
$$
The conclusion of this computation is that
${\partial\over\partial\overline z} Z_3$ is an antiholomorphic 
function of $z$.
Hence $Z_3$ can be written as the sum of two functions $Z_3=h_1+\overline
h_2$ with $h_1$ and $h_2$ both holomorphic.

Consequently if $u$ is a $J$-holomorphic map from $\D$ into
$(\R^6, J)$, $Y_3$ is a harmonic function of $(x,y)$.
If $h$ is a positive harmonic function on $\D$ then we have $|\nabla
h(0)|\le 2h(0)$.  By applying it to $Y_3$, or $-Y_3$, we see that if
$u :\D\to\R^6\setminus M$ is a $J$-holomorphic map $|\nabla
Y_3(0)|\le 2|Y_3(0)|$ (remember $M=\{ y_3=0\}$).  Since $\int^1_0 {dt\over
t} =+\infty$, Proposition 6.1 follows.

\centerline{\bf APPENDICES.}
\bigskip
\noindent {\bf A1. Proof of Proposition 1.1 under the assumption of
$\CC^{k-1,\alpha}$ regularity of $J$.}
\bigskip
Here we simply indicate how to adapt the proof given in section, by using
the trick already used in [N-W]. We had to solve the equation 
${\partial u\over \partial\overline z}+Q_J(u){\partial u\over
\partial z}~=~0$. Let ${\cal B}$ be the closed unit ball in
$\CC^{k,\alpha}(\overline \D, \R^{2n})$.
\smallskip
In order to avoid any differentiation of $J$ in the first step,
we start by fixing a function $\varphi\in {\cal B}$. 
and  we solve instead for $u$ the linear equation:
$${\partial u\over \partial\overline z}+Q_J(\varphi){\partial u\over
\partial z}~=~0\eqno (a.1),$$ 
with prescribed $x$-derivatives  up to order $k$ at 0.
It therefore leads to introducing, in replacement of $\Phi_J$ the operator
$\Phi_J^\#$:
$$\Phi_J^\#(u)~=~({\bf 1}-T_{CG}Q_J(\varphi){\partial \over \partial z})u~,$$
which is simply a linear operator. If
$J$ is close enough to $J_{st}$ in $\CC^{k-1,\alpha}$ topology (independently on $\varphi$)
one can invert this operator. In order to follow the proof of section 1.e
we need to be more precise, we should fix $R>0$ such that 
for any $(q,W)$ in the closed ball ${\bf B}_R$ of radius $R$
in $\R^{2n}\times (\R^{2n} )^k$, $h_{q,W}$ is in the open unit ball of
$\CC^{k,\alpha}(\overline \D,\R^{2n})$. As in section 1.e,
$ h_{q,W}(z)=q+\sum_{l=1}^k{1\over l!}z^lw_l~.$
\smallskip
Then, exactly as in the proof
in section 1.e, to each $(p,V)$ in a neighborhood of 0, independent of
$\varphi \in {\cal B}$,
we associate $(q,W)\in B_R$ such that
$u_{p,V}^\varphi=(\Phi^\#)^{-1} h_{q,W}$ has
the appropriate $k$-jet at 0. Since the function 
$\Phi_J^\#(u_{p,V}^\varphi )$ is an ordinary holomorphic function,
the  function $u=u_{p,V}^\varphi$ satisfies $(a.1)$.
Set $\chi(\varphi )=u_{p,V}^\varphi~.$ If we had
$\varphi=\chi (\varphi )$ then equation $(a.1)$ would give 
${\partial u_{p,V}^\varphi \over \partial\overline z}
+Q_J(u_{p,V}^\varphi){\partial u_{p,V}^\varphi \over
\partial z}~=~0$, so $u_{p,V}^\varphi$ would be $J$-holomorphic
and would solve the problem. 
So we simply need to prove that $\chi$ has a fixed point in
${\cal B}$, provided that $(p,V)$ is close enough to 0.
\smallskip
Note that if $J=J_{st}$,
$Q_J=0$, $\Phi^\#_J$ is the identity, $(p,V)=(q,W)$ and 
$\chi (\varphi)= u_{p,V}^\varphi = h_{p,V}$ (of course independent of $\varphi$).

Since $T_{CG}$ gains one derivative we have the following.
First, if $J$ is close enough to $J_{st}$ in $\CC^{k-1,\alpha}$ topology,
and $(p,V)$ is close enough to 0,
$\chi ({\cal B})\subset {\cal B}$. Second, $\chi$
has a strong continuity property: it maps continuously ${\cal B}$ equipped
with the $\CC^{k-1,\alpha}$ topology into ${\cal B}$ equipped with
the $\CC^{k,\alpha}$ topology. At any rate, $\chi$ defines a continuous map
from ${\cal B}$ equipped with the $\CC^{k-1,\alpha}$ into itself.
Since ${\cal B}$ is a convex compact set in $\CC^{k-1,\alpha}(\overline \D,
\R^{2n})$, The Schauder Fixed Point Theorem implies that $\chi$ has a fixed point,
as desired.
\bigskip\bigskip
\noindent {\bf A2. Deformation of (big) J-holomorphic discs.}
\bigskip
The following Theorem shows the upper semi-continuity
of the Kobayashi Royden pseudo-norm for structures of class 
$\CC^{1,\alpha}$ at least. 
We do not know that it can be reduced
to a question of small perturbation (as the proof of the existence
of $J$-holomorphic discs with prescribed 1-jets for
$\CC^{1,\alpha}$ almost complex structures), and we do not know
a trick similar to the one in  A1 (for dealing with 1-jets and $\CC^\alpha$
almost complex structures) avoiding differentiation of $J$. So
the structure $J$ will be assumed to be of class $\CC^{1,\alpha}$ for 
some $\alpha >0$. Theorem A1 was proved by Kruglikov [K] at least for
structures that are smooth enough. Our proof clarifies the smoothness
assumption and unlike [K] it does not require a careful and difficult 
reading of [N-W].
\bigskip
\noindent {\bf Theorem A1.} 
{\it Let $(X,J)$ be an almost complex manifold with $J$ of
H\" older class $\CC^{1,\alpha}$ ($\alpha >0$).
Let $u$ be a $J$-holomorphic map from a neighborhood of
$\overline \D$ into $X$. There exists a neighborhood $V$ of
$(u(0),{\partial u \over \partial x}(0)$ in the tangent bundle
$TX$ such that for every $(q,Z)\in V$, there exists a
$J$-holomorphic map $v:\D \to X$ with $v(0)=q$ and
${\partial v\over \partial x }(0)=Z.$  }
\bigskip
\noindent {\bf Proof.} Let $r>1$ be such that $u$ is defined on
$\D_r$. 
The map $u$ is $\CC^{2,\alpha}$. 
We can assume that $u$ is an imbedding. Otherwise we add dimensions. We
consider the map $z\mapsto \tilde u (z)= (z,u(z))$ from $\D_r$ into $\R^2\times X$,
equipped with the product almost complex structure $J^2_{st}\times J$.
After getting a map $\tilde v$ from $\D_R$ into $\R^2\times X$, one 
simply takes the projection on the $X$ factor.

By simple
topological arguments, one can find $(n-1)$ smooth vector fields
$Y_1,\ldots ,Y_{n-1}$ defined on a neighborhood of $u(\D_r)$ such that
for every $z\in D_r$ the vectors
${\partial u\over \partial x}(z),Y_1(z),\ldots ,$ $Y_{n-1}(z)$ are
$J(u(z))$-linearly independent.

It allows to define the $\CC^{2,\alpha}$ change of variables
$$(z_1,\ldots ,z_n)\mapsto
u(z_1)+\sum_{j=1}^{n-1}z_{j+1}Y_j(u(z_1))~,$$
defined for $|z|_1<r$ and $|z_j|$ small if $j\geq 2$.

In that change of variables the structure $J$
is transformed into another almost complex structure 
still of class $\CC^{1,\alpha}$, that coincides with the standard one along
$\C\times \{ 0\}\subset \C^n$, The map $u$ is replaced by the map
$z\mapsto (z,0,\ldots ,0)$. So Theorem A1 reduces to the following
Lemma:
\bigskip
\noindent {\bf Lemma A1.} {\it Let $J$ be an almost complex structure on
$\R^{2n}\simeq \C^n$, of H\" older class $\CC^{1,\alpha}$, that coincides with
the standard complex structure on $\C \times \{ 0\}$. Let $U$ be a neighborhood
of $\overline \D \times \{ 0\}$.
For any $(q,t)\in \C^n\times \C^n$ close enough to $(0,0)$, there exists
a $J$-holomorphic map $v:\D \to U$ such that $v(0)=q$ and
${\partial v\over \partial x}(0)=(1,0,\ldots , 0)+t$.}
\bigskip

\noindent {\bf Proof of Lemma A1}. We will work in a neighborhood of
$\overline \D\times \{ 0\}$, on which $J\simeq J_{st}$ and therefore
on which as previously the condition for $J$-holomorphicity can be written:
$${\partial u\over \partial \overline z}
+
Q_J(u) {\partial u\over \partial z}~=~0,$$
and we have $Q_J(z,0,\cdots ,0)=0$.
\bigskip
Assume $J$ of class $\CC^{1,\alpha}$.
Set $\EE_0=\{ f:\overline\D\to \C^n;
f\in\CC^{1,\alpha}, f(0)=0~,~\nabla f(0)=0 \}$,\hfill\break
$\FF_0=\{ g:\overline\D\to \C^n; g\in\CC^\alpha, g(0)=0\}$,
$F(z)=(z,0, \cdots ,0)$.

\noindent Define the map $\Phi$:
$$\eqalign{
\Phi :&\EE_0\longrightarrow \FF_0\cr
    &f\mapsto {\partial (F+f)\over \partial \overline z}
+
Q_J(F+f) {\partial (F+f)\over \partial z}\cr}
$$
Since $F$ is $J$-holomorphic, $\Phi (0)=0$.
We want to show that the derivative of the map $\Phi$
at $f=0$ is onto.
Taking into account that ${\partial F\over \partial\overline z}=0$
and that $Q_J(F)=0$, one gets:
$$\Phi (f) = {\partial f\over \partial \overline z}
+A_J(f)\big( {\partial F\over \partial z}\big)+o(|f|),$$
where $A_J(f)$ is a $(2n\times 2n)$ matrix with entry that are
$\CC^\alpha$ in $z$ and $\R$-linear in $f$.
\bigskip
Denote the derivative at 0 by $D\Phi_0$. With complex instead of real
notations, we can write:
$$
D\Phi_0(f)={\partial f\over\partial\overline z}+B_1(z)f(z)+B_2(z)\overline
f(z)~,
$$
where now $B_1$ and $B_2$ are $(n\times n)$ matrices with complex
coefficients of class $\CC^\alpha$.
\bigskip
The surjectivity of $D\Phi_0$ follows therefore from the following theorem:
\bigskip
\noindent {\bf Theorem A2} {\it If $B_1$ and $B_2$ are $(n\times n)$ complex
matrices with coefficients in $\CC^\alpha (\overline \D)$, for every
$g\in\FF_0$, there exists $f\in \EE_0$ such that
$$
{\partial f\over\partial\overline z}+B_1(z)f(z)+B_2(z)\overline
f(z)=g(z).\eqno(*)
$$}
\medskip

We postpone the proof of Theorem A2, and we now finish the
proof of Lemma A1.  
\bigskip
We shall apply the following elementary result on maps from a Banach space $E$ 
to a Banach space $F$. Let $\Gamma$ be a $\CC^1$ map from $B_E(p,R)$,  the ball
of radius $R$ in $E$ with center at a point $p$, into $F$. Assume that
for some $C>0$,
for all $q\in B_E(p,R)$ the equation $D\Gamma_q(x)=y$ can be solved for
all $y\in F$ with
$\| x\|_E\leq C \| y\|_F$. (By the open mapping theorem, the existence 
of such a constant $C>0$, for $R$ small enough,
is guaranteed as soon as $D\Gamma_0$ is surjective).
Assume moreover that for all $q$ and $q'\in B_E(p,R)$,
$\| D\Gamma_q-D\Gamma_{q'}\|_{\rm op}\leq {1\over 2 C}$. Then for every
$y\in F$, with $\| y-\Gamma (p)\|<  {R\over 2C}$, there exists 
$x\in B_E(p,R)$ such that $\Gamma (x)=y$. The proof is standard, 
by successive approximations.
\bigskip

For $(q,t)$ in $\C^n\times \C^n$, close to $(0,0)$, as in the statement
of the Lemma, and $T=(1,0,\cdots ,0)+t$, set
$$F^\#(z)=q+xT+y\big( J(q)T\big)\qquad (z=x+iy).$$
Then, $F^\#$  close
to $F$ in $\CC^{1,\alpha}$ topology with $\overline\partial_JF^\#(0)=0$.
Define
$$
\Phi^\# : f\mapsto {\partial (F^\# +f)\over \partial \overline z}
+
Q_J(F^\# +f) {\partial (F^\# +f)\over \partial z}
.$$
Since $\overline\partial_JF^\#(0)=0$, $\Phi^\#$ maps
$\EE_0$ into $\FF_0$.
It is a small $\CC^1$ perturbation of $\Phi$ and hence there exists $f$
(close to $0$) such that $ \Phi^\# (f)=0$.
Define $v(z)=F^\# (z)+f(z)$, $v$ has the same first jet at 0 as $F^\#$.
So $\big( v(0),{\partial v\over \partial x}(0)\big)=
(q;(1,\cdots ,0)+t)$ as desired, and
$\Phi^\# (f)=0$ means that  $v$ is $J$-holomorphic. To complete the
proof of Theorem A1, it only remains to prove Theorem A2.
\bigskip
Before starting the proof of Theorem A2, we need 2 Lemmas.
\bigskip

\noindent   {\bf Lemma A2.}  {\it Let $A_1$ and $A_2$ be 
continuous matrices with complex coefficients defined on $\overline\D -\{ 
0\}$ and bounded (no continuity assumed at $0$).  For every continuous map 
$g:\overline\D \rightarrow \C^n$ there exists a continuous map 
$f:\overline D\to \C^n$ such that (in the sense of distribution)
$$
{\partial f\over\partial \overline z} +A_1(z)f+A_2(z)\overline f =g.
$$
}

\noindent (From the proof: We can add $f(0)=0$).

All computations are to be thought in $\R^{2n}$. On  most lines however we 
keep complex $(\C^n)$ notations that allow a simpler writing of 
$\overline\partial$ and of its solution by mean of the Cauchy Kernel.

\medskip

\noindent {\bf Proof of Lemma A2}.  Define the operator $P$:
$$
P(f)={\partial f\over\partial \overline z}+A_1(z)f+A_2(z)\overline f
$$
on $\CC (\overline\D )$ consider the operator:
$$
g\mapsto P\big( g*{1\over\pi z}-g*{1\over \pi z}(0)\big)\;\;\hbox{mapping}
$$
$\CC (\overline\D )$ into itself.  $Q=\11 +K$.  $K$ a compact operator.  
Hence the set of $g\in\CC (\overline D)$ such that $g=P(f)$ for some $f$
is a closed 
subspace ($\R$ subspace not $\C$ subspace) - of finite codimension.

We need to show that it is dense.

Let $\mu$ be a $\R^{2n}$-valued measure on $\overline\D$ such that for 
for any $\varphi\in\CC^\infty (\R^2,\R^{2n})$ such that $\varphi (0)=0$
(the condition  $\varphi (0)=0$ to make sure that $P\varphi$ is continuous at 0):
$$
\int P\varphi\cdot d\mu =0.
$$
We want to prove that $\mu =0$.

By integration by parts $\int P\varphi\cdot  d\mu =0$ gives $P^*\mu =0$ on $\R^2 
-\{ 0\}$ where $P^*$ is an operator of the type:
$$
P^*\mu ={\partial\mu \over\partial z} +C_1(z)\mu +C_2(z)\overline\mu .
$$
The computation is elementary but needs to be done by separating real 
($\R^n$ valued) and imaginary parts of $\mu$.

Since $\mu$ has compact support we first show that $P^*\mu =0$ on 
$\R^2\setminus \{ 0\}$
implies $\mu =0$ on $\R^2\setminus \{ 0\}$.  Note that 
$P^*\mu =0$ on $\R^2 -\{ 0\}$ immediately shows that $\mu\in L^P_{\rm 
loc}$ for any $p<2$ since $\mu ={1\over \pi\overline z}*(-c_1(z)\mu 
+c_2(z)\overline\mu +\nu)$, with $\nu$ a distribution carried by $\{ 0\}$.  
And using the same 
formulas on sees that $\mu$ is a given by a continuous function.  Then
$P^*\mu =0$ gives a 
simple pointwise inequality in $\R^2 -\{ 0\}$
$$
\big|{\partial \mu\over\partial z}\big| <C|\mu |
. $$
By well know uniqueness results it  follows that $\mu=0$ on $\R^2\setminus \{ 0\}$.
\medskip
Consequently, if $\mu$ annihilates the continuous functions which can be written 
as $P(f)$ for $f\in\CC^\alpha$, $\mu$ must be a $2n$-tuple of point masses 
at $0$. Take $f=\overline z a$, $k\in \C^n$, one sees that necessarily $\mu =0$.
\hfill Q.E.D.
\medskip
From now on, the notations are the notations in the statement of
Theorem A2 (for $B_1$, $B_2$ and equation ($*$)).
\medskip
\noindent  {\bf Lemma A3.}  {\it For every $\CC^\alpha$ map 
$g:\overline \D\to \C^n$ with $g(0)=0$, there exists a $\CC^{1,\alpha}$ 
map $f_0$ defined near $0$ in $\C$ such that $f_0(0)=0$, $\nabla f_0(0)=0$ 
and solving ($*$) near $0$.}
\medskip

For $r$ small enough, let $\D_r$ be the disc of radius $r$.  Let $\FF_0(r)$ be the space 
of $\CC^\alpha$ maps $g$ from $\overline \D$ into $\C^n$, with $g(0)=0$.
 To any 
$g\in\FF_0(r)$ associate a $\CC^\alpha$ extension $\widetilde g$ 
with compact support in $\D_{2 r}$ given by a linear operator 
and such that
$$
\sup |\widetilde g|\le 2\sup_{\D_r} |g|.
$$
Now, define $P$ by 
$$
P(f)={\partial f\over\partial\overline z} 
+B_1(z)f+B_2(z)\overline f.
$$
Consider the $\R$-linear map $\Theta$ from $\FF_0(r)$ into itself defined by:

$$\Theta : g\mapsto P\big( \tilde g\ast {1\over \pi z}-(az+b)\big)~,$$
with: $a$ and $b\in \C^n$ are determined by the condition:
$$
\left\{\matrix{\widetilde g*{1\over \pi z}(0)=b\in \C^n\hfill\cr
\cr
{\partial\over \partial z} \big(\widetilde g*{1\over\pi z}\big)(0)=a\in 
\C^n.\hfill\cr}\right.
$$
Due to the choice of $b$, $\Theta$  indeed maps $\FF_0(r)$ into itself. 
It is again a compact perturbation of the identity.  
But  $\Theta$ is clearly one to one if $r$ is small enough.  Indeed, for
$r$ small enough, one has
$$
\sup \big| B_1(z)\big( \widetilde g *{1\over\pi z}-(az+b)\big)\big| +\big| 
B_2(z)\overline{\big(\widetilde g*{1\over\pi z}-(az+b)\big)}\big|\le 
{1\over 2} \sup_{\D_r} |g(z)|.
$$
Therefore $\Theta$ is onto and so one can take $f_0=\widetilde g*{1\over 
\pi z} -(az+b)$. On $\D_r$ $f_0$ solves (*), and $f_0$ is of class
$\CC^{1,\alpha}$. We have $f_0=0$ and due to the choice of $a$,
$\nabla f_0 (0)=0$, as desired.
\medskip\bigskip

\noindent {\bf Proof of Theorem A2.}.  Extend $f_0$ obtained in Lemma A3 , and 
find $f$ by setting $f=f_0+f_1$. 
We need to solve ${\partial f_1\over \partial\overline z}+B_1(z)f_1(z)+B_2(z)\overline 
f_1(z)=g_1(z)$ with $g_1\equiv 0$ near $0$.
It is better to do it on a neighborhood of $\overline \D$ after extending
the data to a neighborhood of $\overline \D$.
Applying Lemma A2 to a larger disc, take $f_1(z)=z^2h(z)$, with \ $h$ obtained by solving
$$
{\partial h\over\partial \overline z}+B_1(z)h(z)+{\overline z^2\over z^2} 
B_2(z)\overline h(z)={g_1(z)\over z^2}.
$$

Then $f=f_0+f_1$ solves ($*$) and satisfies
$$
f(z)=o(|z|)\qquad (z\simeq 0).
$$
Since $f$ solves ($*$), $f$ is of class $\CC^{1,\alpha}$, by elementary 
regularity of $\overline\partial$.

\hfill Q.E.D.

\noindent{\bf A3. Other definition of the Kobayashi distance.}
\smallskip\rm
The usual definition of Kobayashi for the Kobayashi distance
can also be given based on the following
lemma due to Debalme [De-1]
\smallskip\noindent\bf
Proposition A1. {\it Let $J$ be an almost complex structure of class 
$C^{1,\alpha }$ in the neighborhood of the origin in $\R^{2n}$. Then for
any pair $p,q$ of points sufficiently close to the origin there exists a
$J$-holomorphic disk passing through both of them.
}
\smallskip\noindent\bf
Proof. \rm We use again the notations of the proof of Proposition
1.1. But we replace the functions $h_{p,V}$ in the proof of Proposition 1.1
by the functions $h_{p,q}$ dfeined as follows. 

Consider the mapping from $\Delta$ to $\Bbb R ^{2n}$

$$
 h_{p,q}  : \quad  z  \longmapsto p + 2z(q-p)  
$$
where $(p,q) \in (\Bbb R ^{2n})^2$ and denote $u_{\varepsilon ,p,q} 
= \Phi_\varepsilon^{-1} h_{p,q}$. 
We remark that  \hfill \break
$\cdot$ $h_{p,q}$ being holomorphic, $\varepsilon u_{\varepsilon ,p,q} $
is $J$-holomorphic. \hfill \break
$\cdot$ $u_{0,p,q} = h_{p,q}$ . So it verifies $u_{0,p,q}(0)=p$ and
 $u_{0,p,q}({1 \over 2})=q$. \hfill \break
Consider the mapping from $[0, \varepsilon] \times (\Bbb R^{2n})^2$ to $(\Bbb R^{2n})^2$
 $$ \Psi : \quad (\varepsilon,p,q) \longmapsto (u_{\varepsilon,p,q}(0),u_{\varepsilon,p,q}
({1 \over 2}))$$
$\Psi$ is $C^1$ and from our last remark $\Psi(0,.,.) = Id_{(\Bbb R^{2n})^2}$.
So by the implicit functions theorem, if $\varepsilon$ is sufficiently small
, there exists $U$ and $U^{'}$ neighborhoods of zero in $(\Bbb R^{2n})^2$ such that
$ \Psi (\varepsilon,.,.) :  U  \longrightarrow U^{'}$
is a diffeomorphism.
Let $p_0$ and $q_0$ two points sufficiently near of zero (i.e. 
$({p_0 \over \varepsilon},{q_0 \over \varepsilon}) \in U^{'}$). 
There exists $(p,q)$ such that $ \varepsilon u_{\varepsilon,p,q}(0)=p_0$ and 
$\varepsilon u_{\varepsilon,p,q} ({1 \over 2}) = q_0$ . We have thus made 
$ \varepsilon u_{\varepsilon,p,q}$  a $J$-holomorphic curve which is going 
through $p_0$ and $q_0$.
\smallskip\hfill{Q.E.D.}

Both definitions of the Kobayashi distance are equivalent, but we worked only with 
the first one and therefore we will not discuss this matter any further.
\bigskip\noindent
{\bf A4. Classical Properties of the Cauchy-Green Operator.}
\bigskip

For a complex valued function $g$ or a map  $g$ with values in a complex
vector space  continuous on $\overline \D$, and $z\in \C$ with
we set:
$$T_{CG}(g)(z)~=~(g\ast {1\over \pi \zeta})(z)~=~
{1\over \pi} \int_{\D}{g(\zeta)\over z-\zeta}~dxdy(\zeta )~.$$
\bigskip
We have used the following classical properties of $T_{CG}$:
\smallskip\noindent\bf
Proposition A2. {\it For $g$ as abowe:

\noindent (a) ${\partial \over \partial\overline z}[T_{CG}(g)]$ is the distribution 
defined by the function equal to $g$ on $\D$, and to 0 on the complement
of $\D$;

\smallskip\noindent
(b) If $g\in \CC^{k,\alpha}(\overline \D ),~k\in \N~,0<\alpha<1~,$
then the restriction of $T_{CG}g$ to $\overline \D$ belongs to
$\CC^{k+1,\alpha}(\overline \D )$.
}
\bigskip
\noindent Proof. \rm (a) does not need discussion. (b) follows from (a).
The $\CC^{k+1,\alpha}$ regularity of $T_{CG}(g)$ off the unit circle
results from very basic properties of singular integrals
(known as Schauder estimates). See e.g. [M] Chapter II, Theorem 1.6 .
The regularity up to the boundary in each of the regions
$\{ |z|\leq 1\}$ and $\{ |z|\geq 1\}$ is an instance of the
so-called transmission property in the theory of partial
differential equations (Definition 18.2.13 in [H\" o]).
Since it may be harder to find a satisfactory reference
for $\CC^{k+1,\alpha}$ regularity, we  provide a
justification.
\smallskip
Extend $g$ to a function $g_1\in \CC^{k,\alpha}_0(\R^2)$
([St] Chapter VI). Set $h=g_1\ast {1\over\pi z}$.
Then ${\partial h \over \partial \overline z}=g_1$, and
$h\in \CC^{k+1,\alpha}(\R^2)$.
\bigskip
Write $h(e^{i\theta})=h^+(e^{i\theta})+h^-(e^{i\theta})$,
Where $h^+$ is holomorphic on $\{|z|<1\}$, and 
$\CC^{k+1,\alpha}$ on $\{|z|\leq 1\}$; and
 $h^-$ is holomorphic on $\{|z|>1\}$, and 
$\CC^{k+1,\alpha}$ on $\{|z|\geq 1\}$. 
Simply take $h^+$ and $h^-$ to be the Cauchy transform
of the function $e^{i\theta}\mapsto h(e^{i\theta}).$
(Plemelj's formula, and again Schauder's estimates for singular
integrals.)
\bigskip
Consider the function $\tilde h$ defined by:
$$\tilde h = h-h^+~{\rm on}~\overline D,$$
$$\tilde h = h^-~{\rm for}~|z|\geq 1~.$$
Then $\tilde h$ satisfies
$${\partial \tilde h\over \partial \overline z}
={\partial T_{CG}g \over \partial \overline z},$$
on $\R^2$. There is no jump term on the unit circle.
So $\tilde h-T_{CG}g\in \CC^\infty (\R^2)$. Since
$\tilde h$ is by construction $\CC^{k+1,\alpha}$ smooth on
$\overline \D$, so is $T_{CG}g$.\hfill Q.E.D.

\vfill\eject

\centerline{\bf REFERENCES}
\bigskip

\item{[Be]} F.\ Berteloot. Characterization of models in $\C^2$ by their
automorphism groups. Int. J. Math. {\bf 5} (1994) 619-634

\item{[De-1]} R.\ Debalme. Kobayashi hyperbolicity of almost complex manifolds.
Math.CV/98105030 (1998).
\smallskip

\item{[De-2]} R.\ Debalme. Vari\'et\'es hyperboliques presque-complexes.
Th\`ese, Universit\'e de Lille-I (2001).
\smallskip

\item{[D-I]} R.\ Debalme, S.\ Ivashkovich. Complete hyperbolic
neighborhoods in almost complex surfaces.  Int.\ J.\ Math. {\bf 12}
(2001) 211--221.
\smallskip

\item{[Do]} D.\ Donaldson.   Symplectic submanifolds and almost-complex 
geometry.  
J. Differential Geom.  {\bf 44}  (1996),  no. 4, 666--705. 
\smallskip

\item{[G-S]} H.\ Gaussier, A.\ Sukhov. Estimates of the Kobayashi metric
on almost complex manifolds. Math. CV/0307355 (2003).
\smallskip

\item{[Ha]} F.\ Haggui. Fonctions PSH sur une vari\'et\'e presque
complexe. C.R.Acad.Sci.Paris, Ser.I {\bf 335} (2002) 509--514.
\smallskip

\item{[H\" o]} L.\ H\" ormander. {\it The Analysis of Linear Partial
Differential Operators III.} Grund. der math. Wis. 274. Springer-Verlag 
Berlin Heidelberg 1985.
\smallskip

\item{[IS-1]} S.\ Ivashkovich, V.\ Shevchishin.  Structure of the moduli
space in a neighborhood of a cusp curve and meromorphic hulls. 
Invent.\ Math. {\bf 136} (1999) 571--602.
\smallskip

\item{[IS-1]} S.\ Ivashkovich, V.\ Shevchishin. Complex Curves in Almost-Complex
Manifolds and Meromorphic Hulls. ublication Series of 
Graduiertenkollegs "Geometrie und Mathematische Physik" of the Ruhr-University 
Bochum, Issue 36, pp. 1-186 (1999) (see also math.CV/9912046).
\smallskip

\item{[K]} B.S.\ Kruglikov. Existence of Close Pseudoholomorphic
Disks for Almost Complex Manifolds and an Application to the
Kobayashi-Royden Pseudonorm.
Funct. Anal. and  Appl. {\bf 33}
(1999) 38--48.
\smallskip

\item{[K-O]} B.S.\ Kruglikov, M. Overholt. Pseudoholomorphic mappings
and Kobayashi hyperbolicity. Differential Geom. Appl. {\bf 11}
(1999) 265--277.
\smallskip

\item{[Ki]} P.\ Kiernan.  Hyperbolically Imbedded Spaces and Big Picard
Theorem. Math. Ann. {\bf 204} (1973) 203--209.
\smallskip

\item{[M]} S.G.\ Mikhlin. {\it Multidimensional Singular Equations
and Integral Equations.} Pergamon Press (1955). 
\smallskip

\item{[McD]} D.\ McDuff. {Symplectic manifolds with contact type boundaries.
Invent. Math. {\bf 103} (1991) 651-671.
\smallskip

\item{[McD-S]} D.\ McDuff, D.\ Salamon. {\it $J$-holomorphic curves and quantum
cohomology}. Univ.\ Lect.\ Series AMS {\bf 6} (1994).
\smallskip

\item{[N-W]} A.\ Nijenhuis, W.\ Woolf. Some integration problems in
almost complex and complex manifolds. Ann. of Math. {\bf 77} (1963) 424--489.
\smallskip

\item{[Si]} J.-C. Sikorav. Some properties of holomorphic curves in
almost complex manifolds. In {\it Holomorphic Curves in Symplectic
Geometry}, eds. M. Audin and J. Lafontaine, Birkhauser, 1994, 351--361.

\item{[St]} E.M \ Stein {\it Singular Integrals and Differentiability 
Properties of Functions.} Princeton U.P. (1970).

\item{[Za]} M.\ Zaidenberg. Picard's theorem and hyperbolicity. Siberian Math.
J. v.\ 24 (1983) pp. 858-867.

\bigskip
\centerline {- - - - - -  -  - - - - - - - - - - - -}
\bigskip

\nhang{S.\ Ivashkovich: D\'epartement de Math\'ematiques, Universit\'e Lille I,
59655 Villeneuve d' Asq Cedex, France. {\it ivachkov@gat.univ-lille1.fr}}
\smallskip

\nhang{J-P.\ Rosay: Department of Mathematics, University of Wisconsin,
Madison WI 53706 USA. {\it jrosay@math.wisc.edu}}
\bye